\documentclass[reqno]{amsart}
\usepackage[a4paper,left=2.6cm,right=2.4cm,top=2.2cm,bottom=2.5cm]{geometry}
\usepackage[english]{babel}
\usepackage[utf8]{inputenc}
\usepackage[T1]{fontenc}
\usepackage{amsmath, amsthm, amssymb}
\usepackage{graphicx}
\usepackage{subcaption}
\usepackage[hyphens]{url}
\usepackage{booktabs}
\usepackage{lineno}
\usepackage{algorithm}
\usepackage{algpseudocode}

\usepackage{siunitx}
\usepackage{makecell}
\usepackage{colortbl}
\usepackage{xcolor}
\usepackage{tikz}
\usepackage{pgfplots}
\usepackage{pifont}
\usepackage[draft,inline,marginclue]{fixme}
\usepackage{mathtools}
\usepackage[normalem]{ulem}
\usepackage{cancel}
\usepackage{setspace}
\usetikzlibrary{arrows.meta, calc}
\definecolor{forwardblue}{RGB}{39,88,125}
\definecolor{inversegreen}{RGB}{44,108,88}
\definecolor{bluemathlab}{RGB}{0,94,160}
\definecolor{orangemathlab}{RGB}{232,126,26}

\pgfmathdeclarefunction{wig}{2}{%
  \pgfmathparse{1+0.12*sin(3*(#1+#2))+0.06*sin(5*(#1-#2))}%
}

\usepackage [autostyle, english = american]{csquotes}
\MakeOuterQuote{"}
\usepackage[backend=biber,style=numeric,sorting=none,doi=true, isbn=false, url=false, giveninits=true]{biblatex}
\definecolor{darkblue}{rgb}{0.0, 0.0, 0.55}
\definecolor{darkred}{rgb}{0.55, 0.0, 0.0}
\definecolor{checkred}{rgb}{0.75, 0.01, 0.2}
\definecolor{darkgreen}{rgb}{0.0, 0.55, 0.0}

\usepackage{hyperref}
\hypersetup{
    colorlinks=true,
    linkcolor=darkblue,
    filecolor=magenta,
    urlcolor=cyan,
    citecolor=darkblue,
    pdftitle={Inverse domain uncertainty quantification in hemodynamics: sampling virtual patient cohorts under constraints},
    pdfauthor={Francesco Romor, Jose Rodrigo Rojo-Garcia, André-Alexander Zepernick, Ana Djurdjevac, Alfonso Caiazzo},
    pdfkeywords={Bayesian inverse problems; Shape uncertainty quantification; Generative priors; Flow-matching surrogates; Aortic hemodynamics}
}

\newtheorem{theorem}{Theorem}[section]

\theoremstyle{definition}

\theoremstyle{remark}
\newtheorem{remark}[theorem]{Remark}
\theoremstyle{definition}

\newcommand{\KL}{\mathrm{KL}} %

\newcommand{\defeq}{\coloneqq} %
\newcommand{\bx}{\mathbf{x}} %
\newcommand{\bu}{\mathbf{u}} %

\newcommand{\Otheta}{\Omega_{\theta}}
\newcommand{\stheta}{s_{\theta}}
\newcommand{\like}{\Phi_{\text{like}}}
\newcommand{\post}{\Phi_{\text{post}}}
\newcommand{\slike}{\widetilde{\Phi}_{\text{like}}} %
\newcommand{\slikej}[1]{\widetilde{\Phi}_{\text{like},#1}} %
\newcommand{\cdfN}{\Phi_{\mathrm{CDF}}} %
\newcommand{\cdfNinv}{\big(\Phi_{\mathrm{CDF}}\big)^{-1}}

\FXRegisterAuthor{fr}{afr}{FR}
\FXRegisterAuthor{ac}{aac}{AC}
\FXRegisterAuthor{jz}{ajz}{JZ}

\pgfplotsset{compat=1.17}

\begin{document}

\title[Inverse domain UQ]{Inverse domain uncertainty quantification in hemodynamics: sampling virtual patient cohorts under constraints}

\author[Romor, Rojo-Garcia, Zepernick, Djurdjevac, Caiazzo]{Francesco Romor$^{1}$, Jose Rodrigo Rojo-Garcia$^{1}$, André-Alexander Zepernick$^{2}$, Ana Djurdjevac$^{3}$ and Alfonso Caiazzo$^{1}$}

\thanks{$^{1}$Weierstrass Institute for Applied Analysis and Stochastics, Anton-Wilhelm-Amo-Str. 39, 10117 Berlin, Germany.}
\thanks{$^{2}$Institute of Mathematics, Freie Universität Berlin, Arnimallee 6, 14195 Berlin, Germany.}
\thanks{$^{3}$Mathematical Institute, University of Oxford, Woodstock Road, Oxford OX2 6GG, UK}

\begin{abstract}
We consider the following inverse domain uncertainty quantification problem: given a prescribed value or measurement of a hemodynamic biomarker, characterize the distribution of vascular geometries compatible with it. The problem is formulated as a Bayesian inverse problem in which a high-dimensional latent parameter encodes the cardiovascular shape and accounts for inter-patient variability.
The parameter-to-observable map composes a conditional LDDMM stochastic interpolant generating the aortic surface with a high-performance numerical finite element solver for the unsteady incompressible Navier--Stokes equations, from which biomarkers of interest are evaluated. 
The map is nonlinear, high-dimensional and computationally expensive. To address the problem, we compare three ensemble strategies: ensemble adjusted Kalman inversion (EAKI), adaptive Bayesian updating with structural reliability (aBUS), and a Wasserstein--Fisher--Rao (WFR) sequential Monte Carlo scheme. The first two approaches are derivative-free. For the third method, we introduce a novel differentiable flow-matching surrogate that reduces the cost of a single forward map evaluation at the price of a distribution shift depending on the
training distribution used for the surrogate. We apply the framework to the generation of virtual patient cohorts matching prescribed pressure drops along aortic geometries affected by coarctation, and we validate and compare the three methods on synthetic data.

\end{abstract}
\maketitle
\tableofcontents

\section{Introduction}
\label{sec:intro}
Computational blood flow models, in combination with medical imaging, can provide relevant insights into available medical data and support clinical decision making.
Modern numerical methods, together with the constant increase in computational resources and the increasing quality of medical data, contributed to making \textit{in silico} models more and more relevant for several clinical and preclinical applications, including non-invasive diagnosis and staging of cardiovascular pathologies, simulation of pathological conditions and medical devices, as well as treatment planning, optimization, and monitoring. Moreover, the possibility of generating patient-specific and disease-specific synthetic data makes computational models a valuable tool for the development and validation of scientific machine learning algorithms, especially in a field where availability of real data remains a relevant challenge.

However, the quality and reliability of these approaches depend heavily on the model parameters, on the quality of the images used to generate the 
computational domain, and on measurement accuracy.
Uncertainty quantification (UQ) has therefore been increasingly studied for \textit{forward} simulations, i.e., to assess the robustness of hemodynamic predictions with respect to uncertain model inputs, such as boundary conditions, physical, or geometric parameters (see, e.g., ~\cite{sankaran2011stochastic, sankaran2015machinelearning, sankaran2016uncertainty}).
In particular, several recent works focused on studying the impact of geometric uncertainty, such as segmentation errors, on the simulation output,
using machine learning methods to generate geometries from image data~\cite{maher2021convolutionalnetworks}, perturbing a patient-specific aortic surface using Gaussian random fields~\cite{bosnjak2025hemodynamics}, and considering
a generative model, trained on a synthetic patient cohort, 
for generating patient-specific three-dimensional aortic geometries~\cite{KATZ2026119120}.

From the mathematical perspective, the propagation of geometrical uncertainties has been studied more generally in the context of partial differential equations (PDEs) posed on random domains, attracting growing interest in recent years.
In the case of small boundary perturbations, methods from shape calculus and perturbation techniques have been exploited to derive approximations of statistical moments 
of the PDE solution~\cite{harbrecht2013first,harbrecht2008sparse}. 
Alternatively, domain-mapping approaches describe each realization as a deformation of a fixed reference domain and transfer the geometric uncertainty to the coefficients and operators of the transformed PDE~\cite{HPS16,xiu2006numerical}. The latter approach has provided a basis for rigorous regularity and convergence analyses for high-dimensional approximation techniques with respect to the parameters describing the geometry~\cite{djurdjevac2025uncertaintyquantificationstationarytimedependent}.
 
In this work, we consider an inverse domain uncertainty quantification problem, in which observations of the PDE solution are used to infer the underlying domain geometry~\cite{DKOSInverse25, Gantner2018ShapeInversion, Hyvonen2017}. In other words, rather than asking how uncertainty in the vascular geometry affects a given biomarker, the goal is to characterize the distribution of patient geometries compatible with a measured quantity of interest, taking also into account the possible uncertainty in the measurements. 
We formulate this task as a Bayesian inverse problem for uncertain computational domains. 

Focusing on blood flow in the thoracic aorta, we introduce and validate a novel efficient framework for creating virtual patient cohorts with prescribed hemodynamic conditions, to support in silico trials, e.g., for the design and testing of medical devices, and to 
enable cardiovascular shape optimization under geometric uncertainty.

The Bayesian approach tackles both the quantification of uncertainty and the ill-posedness of the inverse problem. For the latter, the regularization of the problem is introduced by modeling the parameter of interest (i.e., the parameters describing the domain) as a random variable distributed according to a \textit{prior} distribution (probability measure). 
Then, in conjunction with a so-called \textit{likelihood function}, which quantifies the error between the results of the forward model and the available measurements, Bayes' theorem is applied to determine a \textit{posterior} distribution of the parameter~\cite{gelman1995bayesian,kaipio2005statistical}.
A rigorous analysis of well-posedness via Bayesian inversion can be found in~\cite{stuart2010inverse,DashtiStuart2017}, where suitable conditions on the mathematical model and prior measure ensure the well-posedness of the posterior measure defined through Bayes' theorem.
The uncertainty quantification process is based on realizations drawn from the posterior measure. From these samples, moments and statistical estimators of the parameter and quantities of interest derived from the forward model can be inferred. Additionally, the samples provide information about the posterior distribution, such as multimodalities and correlations. 
For linear forward models, Gaussian prior distributions and Gaussian additive noise, the posterior distribution is also Gaussian with a closed-form expression for its hyperparameters~\cite{kaipio2005statistical}. In this setting, the sampling process becomes straightforward. 
This is no longer the case for more general nonlinear forward operators, for which classical approaches are based on sampling with Markov chain Monte Carlo methods (MCMC) (e.g.,~\cite{robert1999monte}). However, these techniques may require a high number of evaluations to adequately explore the target distribution, and the cost can be prohibitively high. Several adaptive methods have been proposed to tackle these issues (see, e.g., \cite{christen2010general,haario1999adaptive,haario2001adaptive,haario2006dram}). 

The inverse UQ problem is thus particularly challenging in the present setting. Firstly, one has to handle the high dimensionality of the parameter space (hundreds of degrees of freedom to describe the features of the aortic surface).  Secondly, the map from latent variables to observations is highly nonlinear and involves shape generation, volumetric mesh deformation, and the numerical solution of the unsteady incompressible Navier--Stokes equations. Thirdly, a given observation of a biomarker may be compatible with multiple geometries, potentially leading to multimodal posterior distributions. 

The high-dimensionality of the problem drastically increases the complexity of the sampling process. Two widely used approaches for sampling in high-dimensional spaces are Hamiltonian Monte Carlo (HMC), and Preconditioned Crank-Nicolson (pCN) (see, e.g.,~\cite{cotter2013mcmc}). HMC requires gradient information of the forward model. In contrast, pCN is a derivative-free method, but it may still require a large number of forward model evaluations to adequately explore the posterior distribution. In all cases, the sequential nature of the MCMC methods limits their direct parallelization, and as a consequence, makes simulations within a reasonable time difficult for computationally costly forward models.
A popular alternative to the methods described above is the Laplace approximation, which consists of a local Gaussian approximation around the Maximum a Posteriori point (MAP)~\cite{gelman1995bayesian,helin2022non,wacker2017laplace}. 
By construction, the approximation is expected to perform well for unimodal posterior distributions with a local geometrical structure similar to a Gaussian density. On the other hand, this method requires sufficient differentiability of the likelihood function, which may not always be the case, and it may become computationally expensive in high-dimensional spaces.

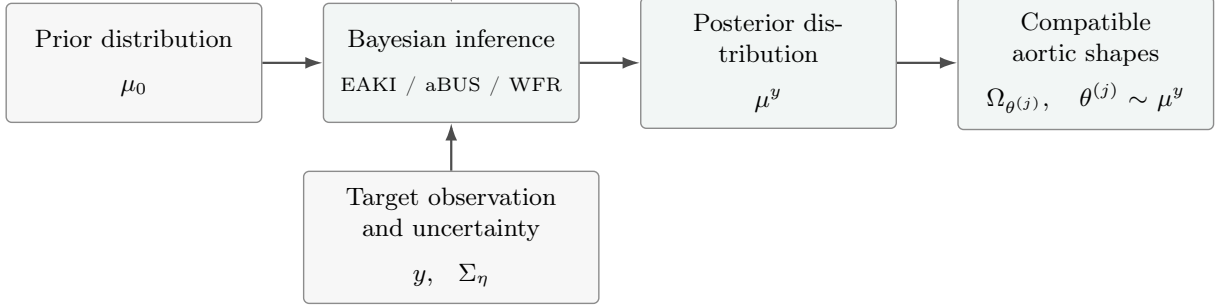
\begin{figure}[htbp]
\centering
\resizebox{\linewidth}{!}{%
\begin{tikzpicture}[
  font=\small,
  box/.style={
    draw=black!45, line width=0.5pt, rounded corners=2pt,
    align=center, text width=2.8cm, minimum height=1.5cm,
    inner sep=6pt, fill=white
  },
  forward/.style={box,fill=forwardblue!6},
  inverse/.style={box,fill=inversegreen!6},
  input/.style={box,fill=black!3},
  arrow/.style={-{Latex[length=2.2mm,width=1.5mm]},
                line width=0.8pt,draw=black!70},
  panel/.style={anchor=west,font=\bfseries\small},
  label/.style={font=\scriptsize,align=center,fill=white,inner sep=2pt}
]

\node[panel,text=forwardblue] at (-1.6,1.35)
  {(a) Forward model $\mathcal G$};

\node[forward] (theta) at (0,0)
  {Latent parameters\\[5pt]$\theta$};
\node[forward] (geometry) at (4,0)
  {Aortic geometry\\[5pt]$\Omega_\theta$};
\node[forward] (flow) at (8,0)
  {Blood flow solution\\[5pt]$(u_\theta,p_\theta)$};
\node[forward] (prediction) at (12,0)
  {Predicted biomarker\\[5pt]$\mathcal G(\theta)$};

\draw[arrow] (theta) -- (geometry);
\draw[arrow] (geometry) -- (flow);
\draw[arrow] (flow) -- (prediction);

\node[label,anchor=north] at (geometry.south)
  {Shape generation \\via LDDMM generative priors};
\node[label,anchor=north] at (flow.south)
  {Unsteady incompressible \\Navier--Stokes equations};
\node[label,anchor=north] at (prediction.south)
  {Pressure-drop curve\\or scalar pressure drop};

\node[panel,text=inversegreen] at (-1.6,-2.45)
  {(b) Bayesian inference};

\node[input] (prior) at (0,-3.6)
  {Prior distribution\\[5pt]$\mu_0$};
\node[inverse] (inference) at (4,-3.6)
  {Bayesian inference\\[5pt]{\scriptsize EAKI / aBUS / WFR}};
\node[inverse] (posterior) at (8,-3.6)
  {Posterior distribution\\[5pt]$\mu^y$};
\node[inverse] (cohort) at (12,-3.6)
  {Compatible aortic shapes\\[5pt]
   $\Omega_{\theta^{(j)}},\quad\theta^{(j)}\sim\mu^y$};

\draw[arrow] (prior) -- (inference);
\draw[arrow] (inference) -- (posterior);
\draw[arrow] (posterior) -- (cohort);

\node[input,text width=3.3cm,minimum height=1.1cm]
  (data) at (4,-5.8)
  {Target observation and uncertainty\\[5pt]$y,\;\Sigma_\eta$};
\draw[arrow] (data.north) -- (inference.south);

\draw[arrow]
  (prediction.south east) -- ++(0,-1.05)
  -| (inference.north);
\end{tikzpicture}}
\caption{Overview of the inverse domain uncertainty quantification framework.}
\label{fig:sketch}
\end{figure}

Our framework builds on two recent works handling arbitrary non-parametric registration of aortic shapes~\cite{romor2026shape} and
generative modeling of geometries based on LDDMM stochastic
interpolants~\cite{KATZ2026119120}, which allow generating
aortic shape distributions from a latent representation of the geometry. Combined with a finite element solver for the
blood flow in the resulting domain, we can then construct a forward map that links the latent geometrical parameter to selected
observations on the fluid solution. In this setting, we investigate three complementary state-of-the-art approaches for high-dimensional Bayesian inversion to solve the inverse domain uncertainty quantification problem (see figure~\ref{fig:sketch}): the Ensemble Adjusted Kalman Inversion (EAKI), the Adaptive Bayesian Updating with Structural reliability (aBUS), and the Wasserstein--Fisher--Rao sequential Monte Carlo (WFR). 

Kalman filters are based on ensembles of particles that iteratively construct Gaussian approximations of the target distribution. The iterative process is derivative-free and typically requires a relatively small number of iterations and model evaluations. 
These methods were initially considered for Bayesian data assimilation of dynamical systems with uncertainty~\cite{sarkka2023bayesian}, and their use has been extended to inverse problems by imposing artificial dynamics, leading to Ensemble Kalman Inversion (EKI) methods (see, e.g.,~\cite{iglesias2013ensemble,iglesias2016regularizing,sanz2023inverse}). 
Inspired by EKI, the EAKI method~\cite{huang2022efficient} introduces suitable artificial dynamics specifically designed for Bayesian inference. Given a Gaussian prior and an ensemble generated from it, the iterative Gaussian densities converge to a Gaussian approximation of the posterior distribution, which, in the case of a linear forward model, coincides with the true posterior density. 

An alternative approach to Bayesian inference is the Bayesian Updating with Structural reliability (BUS) methods~\cite{straub2015bayesian}, which recast the target posterior distribution as a rare event under a standard Gaussian density. This formulation requires rare-event simulation methods such as Subset Simulation (SuS)~\cite{au2001estimation}, which internally relies on MCMC methods. The adaptive BUS (aBUS)~\cite{betz2018bayesian} relies on
an adaptive strategy for the hyperparameters of BUS to improve the speed of convergence.
Unlike the Kalman filter, aBUS does not rely on the assumption of an approximate Gaussian geometry or unimodality of the posterior distribution. Therefore, this approach can efficiently sample complex and multimodal distributions. Moreover, it is derivative-free and, although it internally uses MCMC methods, no burn-in truncation is required. In addition, this method can be efficiently parallelized~\cite{simon2026parallelized}.

Recently, particle-based sampling and optimization methods based on gradient flows have gained significant traction in the context of high-dimensional uncertainty quantification. While classical Langevin dynamics approximate Wasserstein gradient flows to transport particles via diffusion, they often suffer from severe metastability in non-convex or multimodal loss landscapes. Conversely, Fisher–Rao gradient flows enable non-local mass redistribution among different modes. Building on recent studies~\cite{crucinio2025sequential,yan2024learning}, we consider a Wasserstein--Fisher--Rao (WFR) sequential Monte Carlo (SMC) method that combines the two gradient flows (metrics). With respect to MCMC methods (e.g., Hamiltonian Monte Carlo), the proposed approach has a better performance in cases of severe multimodality, such as the one considered in this work (see section~\ref{subsec:results_aortic_wfr}).
Unlike the EAKI and the aBUS, the Wasserstein step of the WFR requires the evaluation of the gradients of the forward map. In the current application, where the forward map combines several complex steps, including
generative modeling for surface deformation,
nonlinear extension of the bulk mesh, and solution of the incompressible Navier--Stokes equations, we circumvent this issue by training a surrogate model that approximates the map and provides cheap batched gradient evaluations on GPUs. 
The machine learning model chosen for this surrogate is an autoregressive flow-matching (stochastic interpolant) model with a transformer-based U-Net, similar to~\cite{bodnar2025aurora}, but extended to the case of variable geometries. To this purpose, we utilize  
an LDDMM registration pre-processing step~\cite{romor2026shape}, which provides meshes with the same topology on different shapes, a natural geometrical encoding that speeds up the training of such heavy architectures.

We test the proposed methods in two different settings: using observations of the pressure drop along the whole vessel centerline and specifying only the maximum value of the pressure drop along the vessel centerline.  The first problem is solved with EAKI, as the observed/specified values of the pressure drop along the vessel centerline are informative enough to reconstruct a single geometry, while aBUS and WFR are used for the second case for sampling multiple geometries associated with the same observed/prescribed value of the maximum pressure drop along the vessel centerline. %
The results and main contributions of this work can be summarized as follows.
\begin{itemize}
  \item We apply three approaches for domain UQ in high-dimensional parameter spaces that do not require gradients of the high-fidelity forward map, comparing them in terms of efficiency and sampling quality.
  \item We show that, even when high-performance numerical solvers are used on relatively coarse meshes, a surrogate forward map provides a remarkable speedup at the price of a small loss of accuracy for the inverse problem under study.
  \item We show that, for topologically equivalent shapes, tailored generative priors and surrogate models that exploit the geometrical encoding provided by registration maps represent a novel and efficient paradigm.
\end{itemize}

The rest of the paper is organized as follows. Section~\ref{sec:inverse_uq} formulates the inverse domain UQ problem and introduces the three considered methods in detail. Section~\ref{sec:forward} describes the forward map which defines the aortic shape inversion problem, as well as the surrogate version used for WFR. Section~\ref{sec:results} presents the numerical results considering a given pressure drop along the aorta centerline and a given maximum pressure drop. Finally, section~\ref{sec:conclusions} draws the conclusions of our study.

\section{Inverse domain uncertainty quantification}\label{sec:inverse_uq}

\subsection{Problem setting}
In what follows, let $\Otheta \subset \mathbb R^d$ ($d=2,3$) denote a computational domain, depending on a parameter vector $\theta \in \mathcal{X} \subset \mathbb{R}^{n_\theta}$, $n_\theta\in\mathbb N$,
where $\mathcal{X}$ stands for a suitable parameter space. 
The relation between $\theta$ and $\Otheta$ 
may correspond to a standard parametrization, such as a free-form deformation or a perturbation induced by a truncated Karhunen--Lo\`eve expansion of a covariance kernel specified on the domain, or, as in the case considered in this work, to a generative model (generative prior)~\cite{xiang2025nativecompactstructuredlatents}, in which $\theta$ plays the role of a latent variable.

On the domain $\Otheta$, we consider a PDE in the general form
\begin{equation}\label{eq:pde}
\mathcal{R}(\Otheta;\stheta) = 0,
\end{equation}
where the solution $\stheta \in V(\Otheta)$ belongs to an abstract function space (e.g., $\stheta=(u_{\theta},p_{\theta})\in H^1(\Otheta)\times L^2(\Otheta)$ for the velocity and pressure fields of the Navier--Stokes equations).
We also define an observation operator $\mathcal O(\theta;\cdot): V(\Otheta) \to \mathbb R^{n_y}$, $\stheta \mapsto \mathcal{O}\bigl(\stheta)$, for $n_y\in\mathbb N$, which extracts
selected quantities of interest from the PDE solution. The observation operator $\mathcal O(\theta;\cdot)$ is assumed to be bounded, not necessarily linear.%
 For the sake of simplicity, we will omit the dependency on $\theta$ in the observation operator in the following description.

The domain uncertainty quantification problem consists in inferring the posterior distribution of $\theta$ given a vector $y \in \mathbb R^{n_y}$ of noisy observations or prescribed values of a quantity of interest. The solution corresponds to a probability distribution on the set of admissible domains $\{\Otheta\}_{\theta\in\mathcal{X}}$. Formally, we introduce a \textit{parameter-to-observable} map
\begin{equation}\label{eq:forward_map}
    \mathcal{G} : \mathcal{X}\subset\mathbb{R}^{n_\theta} \to \mathbb{R}^{n_y}, \qquad \mathcal{G}(\theta) \defeq \mathcal{O}\bigl(\stheta),
\end{equation}
which associates with a particular realization of the parameter $\theta\in\mathcal{X}$ the observation of the corresponding PDE solution $\stheta$ on the domain $\Otheta$.
We model the observations as 
\begin{equation}\label{eq:obs_model}
    y = \mathcal{G}(\theta) + \eta, \qquad \eta \sim \mathcal{N}(0, \Sigma_{\eta}),
\end{equation}
where $\Sigma_{\eta} \in \mathbb{R}^{n_y \times n_y}$ is a symmetric positive-definite noise covariance matrix.

Assuming a Gaussian prior $\theta \sim \mathcal{N}(m_0,C_0)$ on $\mathcal{X}$, with corresponding measure $\mu_0$, Bayes' theorem yields the posterior measure $\mu^y$, which is absolutely continuous with respect to $\mu_0$ and satisfies
\begin{equation}\label{eq:posterior}
    d\mu^{y} = \frac{d\mu^{y}}{d\mu_0}d\mu_{0},
\end{equation}
where
\begin{equation}
\frac{d\mu^{y}}{d\mu_0}(\theta) \propto \exp\!\left(-\frac{1}{2}\bigl\|y - \mathcal{G}(\theta)\bigr\|_{\Sigma_{\eta}^{-1}}^2\right),
\end{equation}
and the weighted norm is defined as $\|\cdot\|_{\Sigma_{\eta}^{-1}}^2 \defeq \langle \cdot, \Sigma_{\eta}^{-1}\cdot\rangle$. 

In the application considered in this work, $\mathcal{G}$ is a nonlinear map acting on a high-dimensional parameter space, which requires the numerical solution of a computationally expensive PDE (see figure~\ref{fig:sketch}). In such situations, the evaluation of the gradient is sometimes computationally prohibitive or undefined if the map is not guaranteed to be differentiable with respect to $\theta$. We therefore focus on derivative-free methods tailored to high-dimensional parameter spaces: an ensemble Kalman inversion method (EAKI) that constructs a sequence of Gaussian approximations %
 (section~\ref{subsec:eaki}), an adaptive Bayesian updating method with structural reliability (aBUS) (section~\ref{subsec:abus}), and a Wasserstein--Fisher--Rao sequential Monte Carlo method (WFR, section~\ref{subsec:wfr}) that leverages a geometry-conditioned LDDMM stochastic interpolant surrogate for efficient gradient evaluation (section~\ref{sec:forward_lddmm}).

For ease of notation, we introduce the potential functions $\like,\post:\mathcal{X}\times \mathbb{R}^{n_y} \to [0,\infty]$, defined as
\begin{equation}\label{eq:philike}
\like(\theta,y):=\frac{1}{2}\bigl\|y - \mathcal{G}(\theta)\bigr\|_{\Sigma_{\eta}^{-1}}^2, 
\end{equation}
and
\begin{equation}\label{eq:phiR}
\post(\theta,y):=\like(\theta,y)+\frac{1}{2}\lVert \theta - m_0 \Vert_{C_0^{-1}}^2.
\end{equation}

\begin{remark}
These methods are not specifically designed for domain UQ problems: for example, they can be employed also when the parameter to be inferred $\theta$ modifies the PDE model, but the domain $\Otheta$ remains fixed. 
Domain UQ problems are characterized by additional challenges, since the parameter-to-observable map requires, besides the resolution of a PDE, additional geometry and mesh processing steps, and a careful design of the shape space taking into account possible geometrical constraints. 
In this work, we consider only geometries with the same topology and for which the
underlying PDE can be defined with the same type of boundary conditions.
Domain UQ and shape optimization problems with varying topologies pose further challenges~\cite{10506652}. 
\end{remark}

\subsection{Ensemble Adjusted Kalman Inversion (EAKI)}
\label{subsec:eaki}
Rather than sampling the posterior $\mu^{y}$ directly, the ensemble Kalman inversion (EKI)~\cite{iglesias2013ensemble, huang2022efficient, chen2024efficient} recasts the inverse problem as a filtering problem for an artificial dynamical system, in which the single available observation $y$ is assimilated repeatedly. The conditional distributions of the state given the data assimilated up to iteration $t$ are approximated by a sequence of Gaussians $\lbrace \mathcal{N}(m_t,C_t) \rbrace_{t \in \mathbb{N}}$, which converges exponentially to the exact posterior mean and covariance if the forward map is linear~\cite[Theorem 1]{huang2022efficient}, and provides a Gaussian approximation of $\mu^{y}$ otherwise, as in our setting.
Only the first two moments of these distributions enter the iteration, and they are estimated empirically from an ensemble of $N$ so-called \textit{particles} ($N\geq 2$).
The method is therefore derivative-free and requires, at each iteration, $N$ independent evaluations of the forward map, which can be run in parallel. At the last iteration $T$ of the artificial dynamics, the empirical mean and covariance of the ensemble provide the estimates $m_{T}$ and $C_{T}$ of the sought Gaussian approximation. Usually, $T$ is specified by looking at convergence metrics, e.g. data misfit or $\Phi_{\text{post}}$.

For the formulation of the method, one considers the augmented operator
\[
\mathcal{F}:\mathcal{X}\to \mathbb{R}^{n_y+n_\theta},
\;\mathcal{F}(\theta) = 
\begin{pmatrix}
\mathcal{G}(\theta) \\
\theta
\end{pmatrix},
\]
and the artificial augmented discrete dynamics
\begin{equation}\label{eq:eaki_dyn}
\left\{
\begin{aligned}
\theta_{t+1} &= \theta_{t} + \varepsilon_{t+1}, \\
z_{t+1} &= \mathcal{F}(\theta_{t+1}) + \nu_{t+1},
\end{aligned}
\right.
\end{equation}
for $t=0,1,\hdots, T$, with an initial Gaussian guess on the parameter $\theta_0 \sim \mathcal{N}(m_0, C_0)$. The second block of the augmented operator $\mathcal{F}$ consists of the parameter itself, so that each iteration assimilates not only the datum $y$, but also the prior mean $m_0$. This definition turns the potential $\like$~\eqref{eq:philike} into the regularized potential $\post$~\eqref{eq:phiR} and prevents the repeated assimilation of the same datum from destroying the prior information.
The noise terms in~\eqref{eq:eaki_dyn}
are sampled from $\nu_{t+1}\sim\mathcal{N}(0, \Sigma_{\nu,t+1})$ and $\varepsilon_{t+1}\sim\mathcal{N}(0, \Sigma_{\varepsilon,t+1})$, respectively, where
\begin{equation}\label{eq:Sigma-nu-t}
        \Sigma_{\nu,t+1}
        :=
        \frac{1}{\Delta\tau}
        \begin{pmatrix}
        \Sigma_{\eta} & 0 \\
        0 & C_0
        \end{pmatrix} 
    \end{equation}
depends on the measurement noise and on the prior on the parameter, and 
\begin{equation}\label{eq:Sigma-eps-t}
        \Sigma_{\varepsilon,t+1}
        :=
        \frac{\Delta\tau}{1-\Delta\tau}\, C_t,
    \end{equation}
depends on the current estimate of the parameter covariance $C_t$. 
The parameter $\Delta\tau\in(0,1)$ in \eqref{eq:Sigma-nu-t} and \eqref{eq:Sigma-eps-t} acts as a time step for the augmented dynamics \eqref{eq:eaki_dyn}: it controls the fraction of the datum and of the prior that is assimilated at each iteration, see also remark~\ref{rmk:fr_flow}.

Let us introduce the augmented observation $z$ defined as 
\begin{equation*}
z_t^{\dagger} := 
\begin{pmatrix}
y \\
m_0
\end{pmatrix},\ \forall t = 0, 1, \ldots%
\end{equation*}
and the observation set at time $t$
\[
Y_t:=\{z_0^{\dagger},\dots,z_{t}^{\dagger}\}.
\]
The method is based on the assumption that the random variables $(\theta, z_{t+1})$ conditioned on $Y_{t}$ are distributed, for each iteration $t$, according to
\begin{equation}
\begin{aligned}
        \theta_t|Y_t&\sim \rho(\theta_t|Y_t) \approx \mathcal{N}(m_t, C_t),\\
        z_{t+1}\mid\theta_{t+1}, Y_t &\sim  \rho(z_{t+1}\mid\theta_{t+1}, Y_t) \propto \rho(\theta_t|Y_t) e^{-\Delta\tau\post(\theta_{t+1},y)},
    \end{aligned}
\end{equation}
where $\post(\theta_{t+1},y)$ is the regularized potential defined in~\eqref{eq:phiR}.

Each iteration is composed of a prediction and an analysis step. In the prediction step, the mean is left unchanged, $\hat{m}_{t+1}=m_t$, while the covariance is inflated to $\hat{C}_{t+1}=C_{t}/(1-\Delta\tau)$, as prescribed by the artificial noise $\varepsilon_{t+1}$ in~\eqref{eq:eaki_dyn}. In the analysis step, the current joint distribution $\rho(\theta_{t+1},z_{t+1} \mid Y_t)$, which is in general not Gaussian, is projected onto the Gaussian distribution that matches its mean and covariance, and a standard Kalman update conditions the latter on the augmented observation $z_t^{\dagger}$. In this way, the conditional distribution $\theta_{t+1}|Y_{t+1}$ is approximated by $\mathcal{N}(m_{t+1}, C_{t+1})$.

Different versions of the method have been proposed, depending on how the ensemble is updated to reproduce the analysis mean and covariance. In the stochastic EKI, each particle is updated using a perturbed observation, and the empirical covariance of the updated ensemble matches $C_{t+1}$ only in a statistical sense, see~\cite[Remark 5]{huang2022efficient}. 
The EAKI employed here belongs instead to the family of square-root filters.
Namely, the particles $\left\{\theta^{j}_{t+1}\right\}_{j=1}^N$ are updated deterministically by applying a linear operator $A:\mathbb{R}^{n_\theta}\rightarrow\mathbb{R}^{n_\theta}$ to the ensemble anomalies, i.e., at each iteration $t$, 
$$\theta^{j}_{t+1} = m_{t+1}+A(\hat\theta^{j}_{t+1}-\hat m_{t+1}),\; j=1,\hdots,N\,.
$$
The operator $A$ is assembled from two singular value decompositions of the ensemble anomalies 
and guarantees that the empirical covariance of the updated ensemble matches $C_{t+1}$ exactly.
Moreover, to prevent the collapse of the particles, we apply a jittering step at the end of each iteration,
in which a perturbation drawn from a centered Gaussian distribution is added to each component of the particles whose standard deviation is below a given threshold.
The EAKI method implemented in this work is described in detail in algorithm~\ref{alg:eaki}.

\begin{remark}[Relation with the Fisher--Rao gradient flow]
    \label{rmk:fr_flow}
    The splitting into a prediction step and an analysis step of the ensemble Kalman inversion method approximates a discrete-time fractional step of the Fisher--Rao gradient flow~\cite{chen2024efficient} minimizing the Kullback--Leibler divergence $\KL(\mu_{t}^{y}\|\mu^{y})$, where $\mu^{y}$ is the posterior distribution defined in~\eqref{eq:posterior} and $\mu_{t}^{y}$ is the one associated with $\theta_t|Y_t$. The fractional step method is approximated as follows:
    \begin{align*}
        \hat{\rho}(\theta_{t+1}|Y_t) &\propto \rho(\theta_t|Y_t)^{1-\Delta\tau}, &\quad&\text{(prediction step)}\\
        \rho(\theta_{t+1}|Y_{t+1}) &\propto \hat{\rho}(\theta_{t+1}|Y_t) e^{-\Delta\tau\post(\theta_{t+1},y)}, &\quad&\text{(analysis step)}
    \end{align*}
    where $\Delta\tau$ is the time step of the fractional step method.
\end{remark}

\begin{algorithm}[htp!]
    \caption{Ensemble adjusted Kalman inversion (EAKI)}
\label{alg:eaki}
\begin{algorithmic}[1]
\setstretch{1.15}
\Require $T_{\text{max}}$ (EAKI iterations), $N$ (ensemble size), $\Delta\tau$ (time step), $\sigma_{\text{thr}}$ (jitter threshold), 
$m_0$, $C_0$ (prior mean and covariance)
\Ensure $m_{T_{\text{max}}}$, $C_{T_{\text{max}}}$ (mean and covariance of the approximate posterior $\theta_{T_{\text{max}}}\sim\mathcal{N}(m_{T_{\text{max}}}, C_{T_{\text{max}}})$).
\State Initialization: sample $\theta_{0}^{j} \sim \mathcal{N}(m_0, C_0)=\mathcal{N}(0, \mathrm{Id}), \; j=1,\hdots,N$
\For{$t = 0,1,\hdots,T_{\text{max}}-1$}
    \Statex \textbf{Prediction step}
    \State $\hat{m}_{t+1} = m_{t}$ 
    \State $\hat{\theta}_{t+1}^{j} = m_{t} + \tfrac{1}{\sqrt{1-\Delta\tau}}(\theta_{t}^{j}-m_{t}), \; j=1,\hdots,N$
    \State $\hat{z}_{t+1}^{j}=\mathcal{F}_{k}(\hat{\theta}_{t+1}^{j}), \; j=1,\hdots,N$ \Comment{parallel forward map}%
    \Statex \textbf{Ensemble generation}
    \State $\hat{z}_{t+1}=\frac{1}{N}\sum_{j=1}^{N}\hat{z}_{t+1}^{j}$
    \State $\hat{Z}_{t+1} = \frac{1}{\sqrt{N - 1}}
\left(
\hat{\theta}_{t+1}^{1} - \hat{m}_{t+1} \quad
\hat{\theta}_{t+1}^{2} - \hat{m}_{t+1} \quad
\cdots \quad
\hat{\theta}_{t+1}^{N} - \hat{m}_{t+1}
\right)$
    \State $\hat{\mathcal{Y}}_{t+1} = \frac{1}{\sqrt{N - 1}}
\left(
\hat{z}_{t+1}^{1} - \hat{z}_{t+1} \quad
\hat{z}_{t+1}^{2} - \hat{z}_{t+1} \quad
\cdots \quad
\hat{z}_{t+1}^{N} - \hat{z}_{t+1}
\right)$
    \State apply SVD to evaluate $P, \hat{D}, V$: $\hat{Z}_{t+1}=P \hat{D}^{1/2}V^{T}$ %
    \State apply SVD to evaluate $U, D$: $V^{T}(\mathrm{Id} + (\hat{\mathcal{Y}}_{t+1})^{T}(\Sigma_{\nu,t+1})^{-1}\hat{\mathcal{Y}}_{t+1})^{-1}V=UDU^{T}$
    \State $A = P\hat{D}^{1/2}UD^{1/2}\hat{D}^{-1/2}P^{T}$
    \Statex \textbf{Analysis step}
    \State $\hat{C}_{t+1}^{\,\theta, z} = \hat{Z}_{t+1}(\hat{\mathcal{Y}}_{t+1})^{T}$
    \State $\hat{C}_{t+1}^{\,z,z} = \hat{\mathcal Y}_{t+1}(\hat{\mathcal Y}_{t+1})^{T} + \Sigma_{\nu,t+1}$
    \State mean: $m_{t+1} = \hat{m}_{t+1} + \hat{C}_{t+1}^{\,\theta,z}\left(\hat{C}_{t+1}^{\,z,z}\right)^{-1}(z_0^{\dagger} - \hat{z}_{t+1})$
    \State variance: $C_{t+1} = \hat{C}_{t+1} - \hat{C}_{t+1}^{\theta z}\left(\hat{C}_{t+1}^{\,z,z}\right)^{-1}\left(\hat{C}_{t+1}^{\,\theta,z}\right)^{T}$  
    \State $\theta_{t+1}^{j} = m_{t+1} + A(\hat{\theta}_{t+1}^{j}-\hat{m}_{t+1})$
    \Statex \textbf{Apply jitter}
    \For{$j= 1, \ldots, N$}
    \For{$l = 1, \ldots, n_\theta$ (components of parameter vector)}
        \If{$\mathrm{std}\bigl((\theta_{t+1})_l\bigr) < \sigma_{\text{thr}}$}
            \State $(\theta_{t+1}^{j})_l = (\theta_{t+1}^{j})_l + \xi, \quad \xi \sim \mathcal{N}(0, 1)$
        \EndIf
    \EndFor
        \EndFor
    \Statex \textbf{Convergence check}
\EndFor
\end{algorithmic}
\end{algorithm}

\subsection{Adaptive Bayesian Updating with Structural reliability (aBUS)}
\label{subsec:abus}

Bayesian Updating with Structural reliability (BUS) methods recast the inference of the posterior
into a rare-event estimation problem, as studied in the context of structural reliability~\cite{straub2015bayesian,betz2018bayesian}. In structural reliability analysis, the state of a system is described by a \textit{limit-state function}, which is negative on the so-called \textit{failure domain}, and the quantity of interest is the failure probability $\mathbb{P}(h\leq 0)$ computed with respect to a standard normal reference density. 

In this context, the difficulty is that failure is typically a rare event. BUS builds on the observation that the posterior can be characterized in the very same form: after augmenting the parameter with one auxiliary variable, the acceptance criterion of rejection sampling defines a domain whose conditional distribution is exactly $\mu^{y}$. 
Unlike EAKI, no Gaussianity assumption is made on the posterior.

To describe the method, let us introduce a bijective map 
$\mathcal{T}: \mathbb{R}^{n_\theta} \to \mathbb{R}^{n_\theta}$ 
that maps the standard normal variables \(u_{1:n_\theta}\) back into the physical parameter space, i.e., such that
$\theta = \mathcal{T}(u_{1:n_\theta})$ for $u_{1:n_\theta}\sim \mathcal{N}(0,\mathrm{Id})$. BUS considers an augmented parameter vector with an independent standard-normal variable $u_{n_\theta+1}\sim\mathcal{N}(0,1)$ and encodes
the observation through the so-called limit-state function~\cite{straub2015bayesian}:
\begin{equation}\label{eq:abus_lsf}
\begin{aligned}
h(u,\ell) & \defeq \log\cdfN(u_{n_\theta+1}) - \log\!\big( e^{-\ell} \,\mathcal{L}(y \mid \mathcal{T}(u_{1:n_\theta}))\big) \\
& = \log\cdfN(u_{n_\theta+1}) + \ell - \log \mathcal{L}(y \mid \mathcal{T}(u_{1:n_\theta})),
\end{aligned}
\end{equation}
where $\mathcal{L}(y \mid \theta)$ is the likelihood function induced by the forward map~\eqref{eq:forward_map}
\begin{equation}\label{eq:abus_lik}
    \mathcal{L}(y \mid \theta) \defeq \exp\!\big(-\like(\theta,y)\big),
\end{equation}
$\ell$ is a scaling constant such that $e^{-\ell} \le 1/\sup_{\theta} \mathcal{L}(y \mid \theta)$, and $\cdfN$ is the standard-normal CDF.
\newcommand{\Omegau}{\mathcal Y_{u}}
Based on~\eqref{eq:abus_lsf}, we define the \textit{failure domain} as the event 
\begin{equation}
    \label{set_ABUS}
    \Omegau := \lbrace u\in\mathbb{R}^{n_\theta+1} \mid h(u,\ell)\leq 0 \rbrace.
\end{equation}
The corresponding conditional probability is given by
\begin{align}\label{eq:abus_cond}
\mathbb{P}\left( u \in \Omegau\mid u_{1:n_\theta}\right) &= \mathbb{P}\left( \cdfN(u_{n_\theta+1}) \le e^{-\ell}\,\mathcal{L}(y \mid \mathcal{T}(u_{1:n_\theta})) \right).
\end{align}
Since $u_{n_\theta+1}\sim \mathcal{N}(0,1)$, we have $\cdfN(u_{n_\theta+1})\sim \mathcal{U}(0,1)$ (standard uniform). As a consequence,  
\begin{equation}
\mathbb{P}(u\in \Omegau \mid u_{1:n_\theta}) =  e^{-\ell} \,\mathcal{L}(y \mid \mathcal{T}(u_{1:n_\theta})),    
\end{equation}
which is equivalent to a standard rejection sampling criterion: a candidate $\theta=\mathcal{T}(u_{1:n_\theta})$ drawn from the prior is accepted whenever the augmented sample $u$ falls in $\Omegau$, that is, with probability $e^{-\ell}\mathcal{L}(y \mid \theta)$. The condition on $\ell$ ensures that this quantity does not exceed one.

Then, we introduce the augmented density given by
\begin{equation}
    \rho_{\Omegau}(u) \propto \mathbf{1}_{\Omegau}(u)\phi_{n_\theta+1}(u),
\end{equation}
where $\phi_{n_\theta+1}$ is the $(n_\theta+1)$-dimensional standard normal density.

Marginalizing the augmented density with respect to $u_{n_\theta+1}$ gives
\begin{equation}
\begin{split}
    \rho_{\Omegau}^{1:n_\theta}(u_{1:n_\theta}) &\propto \phi_{n_\theta}(u_{1:n_\theta})\int_{\mathbb{R}}\mathbf{1}_{\Omegau}(u)\phi_{1}(u_{n_\theta+1})du_{n_\theta+1}\\
    &\propto \phi_{n_\theta}(u_{1:n_\theta})\mathbb{P}(u\in \Omegau \mid u_{1:n_\theta})\\
    &\propto \phi_{n_\theta}(u_{1:n_\theta})\mathcal{L}(y \mid \mathcal{T}(u_{1:n_\theta})).
\end{split}
\end{equation}
The marginal density thus corresponds to the sought density of $\mu^{y}$ under a change of coordinates.
Hence, to obtain samples of $\mu^{y}$, the first $n_\theta$ components of the samples drawn from the augmented density are pushed through the map $\mathcal{T}$.
This formulation in the augmented space is particularly convenient, since the resulting density is expressed with respect to a standard normal reference density, while the prior information is preserved through the map $\mathcal{T}$.

Sampling $\mu^{y}$ by plain rejection requires, on average, one accepted sample every $\big(e^{-\ell}\int_{\mathcal{X}} \mathcal{L}(y \mid \theta)\,d\mu_{0}\big)^{-1}$ evaluations of the forward map.
For informative data, the likelihood is concentrated with respect to the prior, the domain $\Omegau$ becomes a rare event under the standard normal reference, and this cost, like that of standard MCMC methods in high-dimensional spaces, becomes prohibitive.
We therefore consider the combination of BUS
with Subset Simulation (SuS)~\cite{au2001estimation} proposed in~\cite{straub2015bayesian}, in which the rare event is reached through a sequence of nested and progressively less probable events.

Namely, we consider a series of intermediate threshold levels
\[
+\infty=: h_0>h_1>\dots>h_{n_L}:=0,
\]
and define the nested events 
\[
F_{i} := \lbrace u \mid h(u,\ell)\leq h_{i} \rbrace,\;i=0,1,\hdots,n_L\,,
\]
as well as the densities
\begin{equation}
    \rho_{i}(u) \propto \mathbf{1}_{F_{i}}(u)\phi_{n_\theta+1}(u),\;i=0,1,\hdots,n_L.
\end{equation}
By definition, the level sets $F_i$ form a nested chain, i.e., 
\[
\mathbb R^{n_\theta+1} = F_0 \supseteq F_1 \supseteq \cdots \supseteq F_{n_L} = \Omegau.
\]
Hence, the acceptance probability telescopes into the product of the level-conditional probabilities, i.e.,
\begin{equation}\label{eq:abus_evidence}
    \mathbb{P}(h\le 0)=\mathbb{P}(F_{n_L})=\prod_{k=1}^{n_L}\mathbb{P}(F_k\mid F_{k-1})=\prod_{k=1}^{n_L}P_k.
\end{equation}

In practice, at each level $i = 0,1,\hdots,n_L-1$, we adaptively select $h_{i+1}$ as the $p_{0}$-quantile of the ordered samples at level $i$, so that a fraction $p_0$ of them satisfies the event $F_{i+1}$. These samples are retained as seeds and used to initialize multiple MCMC chains with target distribution $\rho_{i+1}$.
By construction, the seeds already belong to $F_{i+1}$ and are distributed according to $\rho_{i+1}$: the chains therefore start in stationarity and no burn-in truncation period is required.

The efficiency of the MCMC sampling process depends on the proposal distribution. 
In the context of SuS, some specialized methods have been proposed (see, e.g.,~\cite{botev2012efficient,cerou2012sequential,del2006sequential}), 
in particular, we use adaptive conditional sampling (aCS) as introduced in~\cite{papaioannou2015mcmc}.

The overall efficiency of BUS in conjunction with SuS depends on the choice of the proportionality
constant $e^{-\ell}$: the smaller it is, the rarer the event $\Omegau$ and the larger the number of levels, while its admissibility condition $e^{-\ell}\le 1/\sup_{\theta}\mathcal{L}(y\mid\theta)$ involves the maximum of the likelihood, which is not known a priori.
The adaptive BUS (aBUS) method addresses this issue by estimating the constant at runtime~\cite{betz2018bayesian}: $\ell$ is replaced by the largest log-likelihood value $\hat\ell$ observed so far and, whenever a larger value is found, the current threshold is shifted accordingly and the auxiliary variable $u_{n_\theta+1}$ is resampled from its conditional distribution (last two steps of the loop in algorithm~\ref{alg:abus}), so that the samples remain distributed according to the current intermediate density.

Each evaluation of the limit-state function~\eqref{eq:abus_lsf} requires one evaluation of the forward map $\mathcal{G}$, which includes the numerical solution of the Navier--Stokes equations on the different domains (section~\ref{sec:forward}).
Since the candidates of a level, or of a batch of MCMC steps, are mutually independent, they can be evaluated concurrently: following~\cite{simon2026parallelized}, each batch of forward evaluations $\theta\mapsto\mathcal{G}(\theta)$ is dispatched at once, so that the $N$ solutions of the PDE underlying the forward map are run in parallel.
The implemented aBUS algorithm is summarized in algorithm~\ref{alg:abus}.%

\begin{algorithm}[t]
\caption{Parallel adaptive BUS with subset simulation (aBUS--SuS). The forward map is detailed in 
section~\ref{sec:forward} and algorithm~\ref{alg:g-theta}.}
\label{alg:abus}
\begin{algorithmic}[1]
\setstretch{1.15}
\Require samples per level $N$, conditional probability $p_0=P(F_i|F_{i-1})$, with $i=1, \dots,n_L$, prior transform
  $\mathcal{T}(\tilde\theta)=m_c+L_c\tilde\theta$, log-likelihood $\log \mathcal{L}(y \mid \theta)=-\like(\theta,y)$
\Ensure thresholds $\{h_i\}$, posterior samples $\{\theta_j\}$
\State $u=(u_{1:n_\theta},u_{n_\theta+1})\in\mathbb{R}^{n_\theta+1}$, with $u_{n_\theta+1}\sim\mathcal{N}(0,1)$ \Comment{augment the standard space}
\State $h(u,\hat\ell) = \log\cdfN(u_{n_\theta+1}) + \hat\ell - \log \mathcal{L}\big(y \mid \mathcal{T}(u_{1:n_\theta})\big)$ \Comment{limit-state function~\eqref{eq:abus_lsf}}
\Statex
\Statex \textbf{Initial Monte Carlo level}
\State draw $u_j \sim \mathcal N(0,\mathrm{Id})$, \ $j=1,\dots,N$
\State $\{\log \mathcal{L}_j\} = \Call{ForwardMap}{\{u_j\}}$;  \Comment{forward map, in parallel}
\State $\hat\ell = \max_j \log \mathcal{L}_j$; \ \ $i = 0$; \ \ $h_0 = +\infty$
\Statex
\Statex \textbf{Subset simulation}
\While{$h_i > 0$}
  \State $i = i+1$; \ \ $g_j = h(u_j,\hat\ell)$ for all $j=1,\dots,N$
  \State $h_i =$ $p_0$-quantile of $\{g_j\}$; \ \ $n_F = \#\{\,j : g_j \le \max(h_i,0)\,\}\approx p_0\cdot N$
  \State \textbf{if} $h_i < 0$ \textbf{then} $h_i = 0$, \ $P_i = n_F/N$ \ \textbf{else} $P_i = p_0$
  \State $\{\omega_k\}_{k=1}^{n_F}$ $=$ the $n_F$ samples with smallest $g_j$ \Comment{set up the seeds}
  \State $\big(\{u_j\},\{\log \mathcal{L}_j\},\lambda\big) = \Call{aCS}{N,\lambda,h_i,\{\omega_k\}_{k=1}^{n_F},\hat\ell}$ \Comment{adaptive MCMC, algorithm~\ref{alg:abus-sub}}
  \State $\ell' = \max(\hat\ell,\ \max_j \log \mathcal{L}_j)$; \ \ $h_i = h_i - \hat\ell + \ell'$; \ \ $\hat\ell = \ell'$ \Comment{update threshold}
  \State $u_{j,n_\theta+1} = \cdfNinv(p_j)$, \ \ $p_j \sim \mathcal U\!\big(0,\ \min(1,\ e^{\log \mathcal{L}_j - \hat\ell + h_i})\big)$ \Comment{refresh BUS variable}
\EndWhile
\State $n_L = i$; \ \ $\theta_j = \mathcal{T}(u_{j,1:n_\theta})$
\end{algorithmic}
\end{algorithm}

\begin{algorithm}[t]
\caption{Adaptive conditional sampling (aCS) called by algorithm~\ref{alg:abus}.
The forward map is detailed in 
section~\ref{sec:forward} and algorithm~\ref{alg:g-theta}.}
\label{alg:abus-sub}
\begin{algorithmic}[1]
  \Require target acceptance value ($\hat{a}_{\text{target}}=0.44$), initial proposal std ($\sigma_{\mathrm{aCS},0}=1$), batch size ($n_B=100$).
\Function{aCS}{$N,\ \text{proposal scaling parameter } \lambda,\ \tau,\ \{\omega_k\}_{k=1}^{n_F},\ \hat\ell$}
  \State distribute $N$ samples over $n_F$ chains; \ seed each chain head from $\{\omega_k\}_{k=1}^{n_F}$
  \State $\sigma_{\mathrm{aCS}} = \min(\lambda\sigma_{\mathrm{aCS},0},\mathbf 1)$; \ \ $\varrho = \sqrt{1-\sigma_{\mathrm{aCS}}^2}$
  \State split the $N$ chain steps into $n_A=\lceil N/n_B\rceil$ batches of (at most) $n_B$ proposals each
  \For{$k = 1,\dots,n_A$}
    \State propose $v^\star \sim \mathcal N\!\big(\varrho\,u_{t-1},\ \sigma_{\mathrm{aCS}}^2\big)$ for every active chain \Comment{conditional proposal}
    \State $\{\log \mathcal{L}^\star\} = \Call{ForwardMap}{\{v^\star\}}$ \Comment{forward map, in parallel}
    \State accept $u_t = v^\star$ if $h(v^\star,\hat\ell) \le \tau$, else $u_t = u_{t-1}$ (carry $\log \mathcal{L}$)
    \State evaluate batch acceptance rate $\hat a$
    \State $\lambda = \exp\!\big(\log\lambda + \tfrac{1}{\sqrt{k}}(\hat a - \hat{a}_{\text{target}})\big)$; \ update $\sigma_{\mathrm{aCS}},\varrho$ 
  \EndFor
  \State \Return $\{u_j\}$, $\{\log \mathcal{L}_j\}$, $\lambda$
\EndFunction
\end{algorithmic}
\end{algorithm}

\subsection{Wasserstein--Fisher--Rao sequential Monte Carlo (WFR)}
\label{subsec:wfr}

As pointed out in remark~\ref{rmk:fr_flow}, EAKI approximates a fractional step scheme of the Fisher--Rao gradient flow of the
Kullback--Leibler divergence $\KL(\mu_t^y\|\mu^y)$. 
In this section, we introduce the Wasserstein--Fisher--Rao sequential Monte Carlo (WFR) method, which
can be seen as a further extension~\cite{yan2024learning,zhu2026inclusiveklgradientflows} of EAKI or of the Fisher--Rao gradient flow (see remark~\ref{rmk:fr_flow})
that approximates the full Wasserstein--Fisher--Rao gradient flow of $\KL(\rho\|\mu^y)$, 
for a probability measure $\rho$.

The idea behind the method is to evolve an empirical measure $\rho_t$ towards the posterior $\mu^y$ of~\eqref{eq:posterior} along the gradient flow of the Kullback--Leibler divergence $\KL(\rho\|\mu^y)$ under the Wasserstein--Fisher--Rao metric~\cite{yan2024learning,zhu2026inclusiveklgradientflows},
combining the Wasserstein (optimal-transport) and Fisher--Rao (birth--death) geometries, i.e.,
\begin{equation}\label{eq:wfr_flow}
    \partial_t\rho_t = \underbrace{\nabla\!\cdot\!\big(\rho_t\,\nabla\log(\rho_t/\mu^y)\big)}_{\text{Wasserstein transport}} \;-\; \underbrace{\rho_t\Big(\log\tfrac{\rho_t}{\mu^y}-\mathbb E_{\rho_t}\!\big[\log\tfrac{\rho_t}{\mu^y}\big]\Big)}_{\text{Fisher--Rao birth--death}}.
\end{equation}
The measure 
is approximated by an ensemble $\{\theta_j\}_{j=1}^N$ of $N$ particles, i.e., $\rho_t\approx\frac1N\sum_{j}\delta_{\theta_j}$, and the evolution is realized
via an operator splitting of the two half-flows over a step $\gamma_t$~\cite{crucinio2025sequential}. 
Unlike EAKI and aBUS, WFR requires the evaluation of the likelihood score, i.e., the gradient of the log-likelihood with respect to the parameters, which in turn requires the computation of the Jacobian $\partial\mathcal{G}/\partial\theta$ of the forward map.
To reduce the computational cost, we replace $\mathcal{G}$ with a differentiable neural network surrogate $\Tilde{\mathcal G}$ (section~\ref{sec:forward_lddmm}) that approximates the PDE solution on the deforming domain.

Let 
$$
\slike(\theta, y):=\tfrac12\lVert y-\Tilde{\mathcal G}(\theta)\rVert_{\Sigma_\eta^{-1}}^2
$$
denote the potential~\eqref{eq:philike} evaluated through the differentiable surrogate $\Tilde{\mathcal G}$. Writing the log-posterior, up to an additive constant, as $\log\mu^y(\theta)=-\slike(\theta, y)+\log\pi_0(\theta)$, 
and the prior $\pi_0=\mathcal N(0,\mathrm{Id})$, and approximating the ensemble density $\rho_t$ by the Gaussian-mixture kernel estimate $q(\theta)=\frac1N\sum_i\mathcal N(\theta;\theta_i,2\gamma_t \mathrm{Id})$, where $\mathcal N(\theta;\theta_i,2\gamma_t \mathrm{Id})$ denotes the density of the Gaussian distribution with mean $\theta_i$ and covariance $2\gamma_t \mathrm{Id}$ evaluated at $\theta$, one analysis step applies the following updates in sequence:
\begin{itemize}
    \item \emph{Fisher--Rao (birth--death).} Each particle is scored by the log density-ratio 
    $$
    \ell_j=\log\mu^y(\theta_j)-\log q_j=-\slikej{j}+\log\pi_0(\theta_j)-\log q_j
    $$
    between the posterior and the ensemble, positive where $\theta_j$ is under-represented relative to $\mu^y$ and negative where it is over-represented. The particles are reweighted and normalized,
    \begin{equation}\label{eq:wfr_weights}
        w_j=\frac{\exp(\alpha_{\mathrm{FR}}\,\ell_j)}{\sum_{k=1}^N\exp(\alpha_{\mathrm{FR}}\,\ell_k)},
    \end{equation}
    then resampled (systematic resampling), duplicating under-represented particles and pruning over-represented ones. The exponent $\alpha_{\mathrm{FR}}$ is the Fisher--Rao step size: its nominal value $1-e^{-\gamma_t}$ discretizes the birth--death half-flow over $\gamma_t$, but to prevent weight degeneracy it is decreased to the largest $\alpha_{\mathrm{FR}}\le 1-e^{-\gamma_t}$ %
    for which the effective sample size
    \begin{equation}\label{eq:wfr_ess}
        \mathrm{ESS}(\alpha_{\mathrm{FR}})=\Big(\textstyle\sum_{j=1}^N w_j^2\Big)^{-1}\in[1,N],
    \end{equation}
    stays above $N/2$; $\mathrm{ESS}=N$ for uniform weights and $\mathrm{ESS}=1$ when a single particle dominates.
    \item \emph{Wasserstein (Langevin).} Each resampled particle is transported by one Langevin step
    \begin{equation}\label{eq:wfr_langevin}
        \theta_{t+1}=\theta_t+\gamma_t\big({-\nabla_\theta\slike}+\nabla_\theta\log\pi_0\big)+\sqrt{2\gamma_t}\,\xi_t,\qquad \xi_t\sim\mathcal N(0,\mathrm{Id}),
    \end{equation}
    with likelihood score 
    $$
    -\nabla_\theta\slike(\theta)=\big(\partial\Tilde{\mathcal G}/\partial\theta\big)^{\!\top}\Sigma_\eta^{-1}\big(y-\Tilde{\mathcal G}(\theta)\big), 
    $$
    evaluated from the Jacobian $\partial\Tilde{\mathcal G}/\partial\theta$ of the surrogate by automatic differentiation (see remark~\ref{rmk:grad}), and prior score $\nabla_\theta\log\pi_0=-\theta$, so that the drift equals $\nabla_\theta\log\mu^y$ and the invariant law of the continuous dynamics is the posterior $\mu^y$.
\end{itemize}

Unlike the ensembles of sections~\ref{subsec:eaki}--\ref{subsec:abus}, the approximate forward map $\Tilde{\mathcal G}$ and its Jacobian $\partial\Tilde{\mathcal G}/\partial\theta$ are here evaluated in a single batched, differentiable forward pass of the PDE surrogate introduced in section~\ref{sec:forward_lddmm}.
Consequently, the flow does not require a separate high-fidelity PDE solution 
for each ensemble member. An additional advantage of a neural network surrogate is the use of batched evaluations of observables and their gradients, which is more efficient than the sequential evaluation of the high-fidelity PDE solver. Algorithms~\ref{alg:wfr} and~\ref{alg:wfr-sub} sketch the
steps of the WFR method.

\begin{algorithm}[t]
\caption{Parallel Wasserstein--Fisher--Rao (WFR) sequential Monte Carlo sampling. The surrogate forward map is defined in section~\ref{sec:forward_lddmm}, see also algorithm~\ref{alg:g-tilde-theta}. The analysis step is detailed in algorithm~\ref{alg:wfr-sub}.}
\label{alg:wfr}
\begin{algorithmic}[1]
\Require ensemble size $N$, iterations $T_{\max}$, noise scale $\sigma_\eta$ ($\sigma_\eta^2=\Sigma_{\eta}$ for the scalar observable $\mathcal{O}_2$),
  observed datum $y$, parameter space maps $(L_c,m_c,U)$ of~\eqref{eq:latent_reconstruction},
  differentiable surrogate $\Tilde{\mathcal G}$ of $\mathcal G$
\Ensure posterior ensemble $\{\Tilde\theta_j\}_{j=1}^N$ 
\State $\Tilde{\theta}_j \sim \mathcal N(0,\mathrm{Id})$ \Comment{prior init $\pi_0=\mathcal N(0,\mathrm{Id})$}
\For{$t = 0,\dots,T_{\max}-1$}
  \State $\gamma_t = \gamma_{\min} + \tfrac12(\gamma_{\max}-\gamma_{\min})\big(1+\cos(\pi\, t/(T_{\max}-1))\big)$ \Comment{cosine step scheduler}
  \State $\{\Tilde{\mathcal G}(\Tilde{\theta}_j)\},\{\nabla\Tilde{\mathcal G}(\Tilde{\theta}_j)\} = \Call{SurrogateForwardMap}{\{\Tilde{\theta}_j\}}$ \Comment{batched differentiable GPU job}
  \State $\{\Tilde{\theta}_j\} = \Call{WFRanalysis}{\{\Tilde{\theta}_j\},\{\Tilde{\mathcal G}(\Tilde{\theta}_j)\},\{\nabla\Tilde{\mathcal G}(\Tilde{\theta}_j)\},\,y,\sigma_\eta,\gamma_t}$
\EndFor
\State \Return $\{\Tilde{\theta}_j\}$
\end{algorithmic}
\end{algorithm}

\begin{algorithm}[t]
\caption{One WFR sequential MC step called by
algorithm~\ref{alg:wfr}. 
}
\label{alg:wfr-sub}
\begin{algorithmic}[1]
\Function{WFRanalysis}{$\{\Tilde{\theta}_j\},\{\Tilde{\mathcal G}(\Tilde{\theta}_j)\},\{\nabla\Tilde{\mathcal G}(\Tilde{\theta}_j)\},\,y,\sigma_\eta,\gamma$}
\Statex \hspace{-1.2em}\textit{(1) Likelihood potential and its score}
  \State $r_j = \tilde{\mathcal G}(\tilde\theta_j) - y$; $\slikej{j} = \tfrac12 \big(r_j/\sigma_\eta\big)^2$
    \Comment{Gaussian likelihood potential~\eqref{eq:philike}}
  \State $\mathbf s_j = -\nabla_{\tilde\theta}\slikej{j}
          \;=\; -\,\frac{r_j}{\sigma_\eta^{2}}\,\nabla_{\tilde\theta}\tilde{\mathcal G}(\tilde\theta_j)$
    \Comment{likelihood score}

\Statex \hspace{-1.2em}\textit{(2) Fisher--Rao flow: birth--death by reweighting and resampling}
  \Statex \hspace{1.2em} the FR velocity is $\log(\mu^y/\rho)$, with target
    $\mu^y \propto \pi_0\,e^{-\slike}$ and $\rho$ the current particle density
  \State $\log q_j = -\log N +
          \operatorname*{logsumexp}_{1\le i\le N}
          \Big[-\tfrac{n_\theta}{2}\log(4\pi\gamma) - \tfrac{1}{4\gamma}\lVert \tilde\theta_j - \tilde\theta_i\rVert^2\Big]$
    \Comment{$q \approx \rho$}
    \Statex \hspace{1.2em} a Gaussian-mixture density estimate over the ensemble, kernel $\mathcal N(\cdot\,;\tilde\theta_i, 2\gamma\mathrm{Id})$
  \State $\log \pi_0(\tilde\theta_j) = -\tfrac12 \lVert \tilde\theta_j\rVert^2$
    \Comment{$\tilde\theta \sim \mathcal N(0,\mathrm{Id})$ prior, up to an additive constant}
  \State $\ell_j = -\slikej{j} + \log\pi_0(\tilde\theta_j) - \log q_j$
    \Comment{$\log\mu^y - \log q$; constants cancel in $w$ below}
  \State $\bar w_j(\alpha_{\mathrm{FR}}) =
          \exp(\alpha_{\mathrm{FR}} \ell_j)\big/\textstyle\sum_{i} \exp(\alpha_{\mathrm{FR}} \ell_i)$, \ \
          $\operatorname{ESS}(\alpha_{\mathrm{FR}}) = \big(\textstyle\sum_j \bar w_j(\alpha_{\mathrm{FR}})^2\big)^{-1} \in [1,N]$
  \If{$\operatorname{ESS}(1) \ge N/2$}
    \State $\alpha_{\mathrm{FR}} = 1$
  \Else
    \State $\alpha_{\mathrm{FR}} =$ largest $\alpha_{\mathrm{FR}} \in [0,1]$ with $\operatorname{ESS}(\alpha_{\mathrm{FR}}) \ge N/2$, by $40$ bisection steps
      \Comment{$\operatorname{ESS}$ decreases in $\alpha_{\mathrm{FR}}$}
  \EndIf
  \State $w_j = \bar w_j(\alpha_{\mathrm{FR}})$
    \Comment{$w_j \propto (\mu^y/q)^{\alpha_{\mathrm{FR}}}$ at $\tilde\theta_j$}
  \State draw $u \sim \mathcal U(0,1)$; \ \
         $\operatorname{idx}_j = \min\big\{\,i:\ \sum_{l \le i} w_l \ \ge\ (j+u)/N \,\big\}$, \ $j=0,\dots,N-1$
    \Comment{systematic resampling}
  \State $\tilde\theta_j = \tilde\theta_{\operatorname{idx}_j}$, \ \
         $\mathbf s_j = \mathbf s_{\operatorname{idx}_j}$
    \Comment{scores travel with their particles}

\Statex \hspace{-1.2em}\textit{(3) Wasserstein flow: one Langevin step on the resampled particles}
  \State $v_j = \mathbf s_j - \tilde\theta_j$
    \Comment{$\nabla_{\tilde\theta}\log\mu^y$: likelihood score $+$ prior score $\nabla_{\tilde\theta}\log\pi_0 = -\tilde\theta$}
  \State $\xi_j \sim \mathcal N(0, \mathrm{Id})$; \ \
         $\tilde\theta_j = \tilde\theta_j + \gamma\, v_j + \sqrt{2\gamma}\;\xi_j$
    \Comment{Euler--Maruyama, unadjusted}
  \State \Return $\{\tilde\theta_j\}_{j=1}^{N}$
\EndFunction
\end{algorithmic}
\end{algorithm}

\section{Forward problem}
\label{sec:forward}

\begin{figure}[!ht]
  \includegraphics[width=0.9\textwidth]{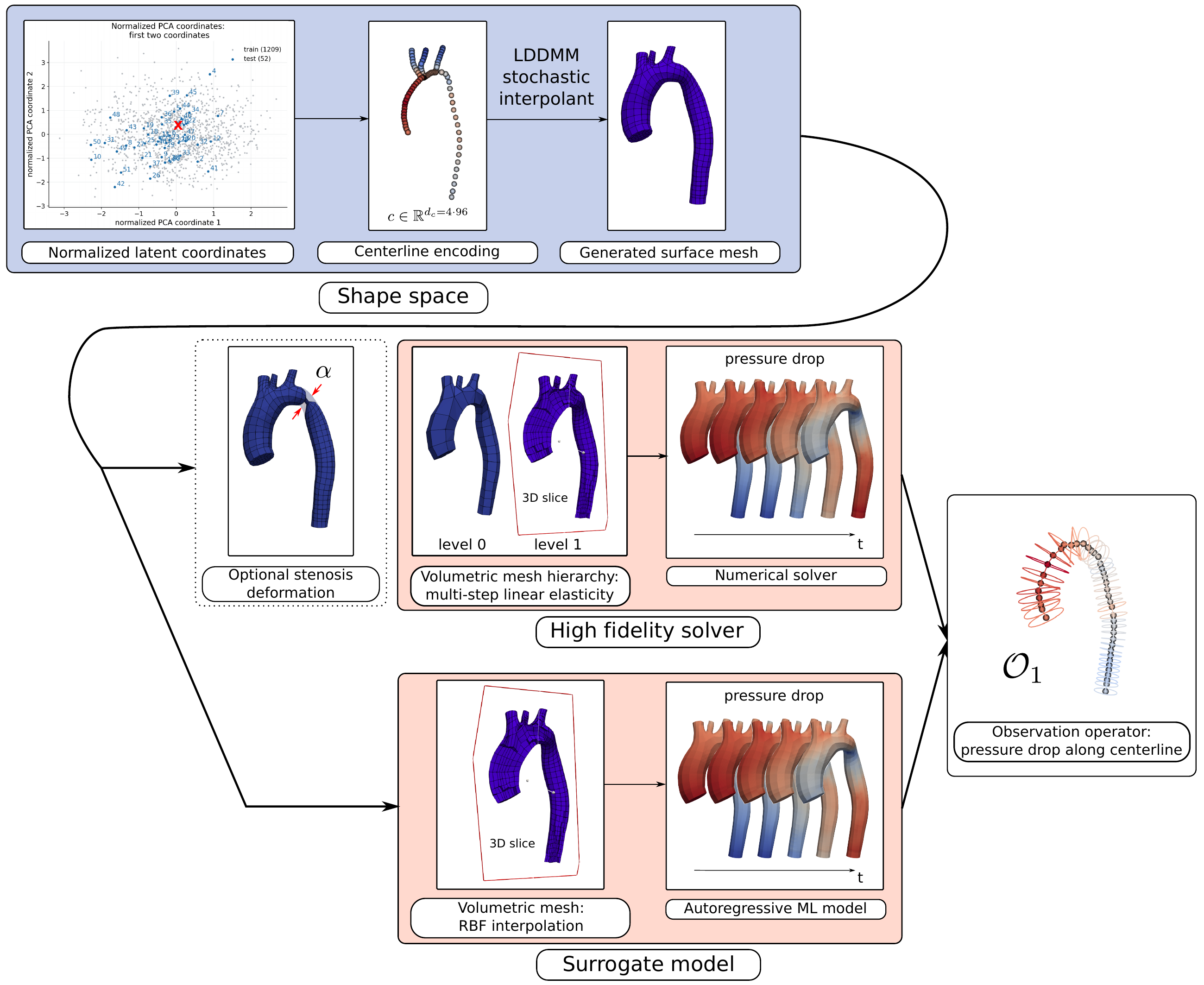}
  \caption{Schematic representation of the forward map, from the latent representation of the shape to the observation, using either a high-fidelity numerical solver (section~\ref{ssec:map}) or a surrogate model (section~\ref{sec:forward_lddmm}).}
  \label{fig:scheme_forward}
\end{figure}

This work is motivated by the problem of finding a posterior distribution compatible with a measured hemodynamic biomarker in aortic blood flows, such as pressure drops.
In the formalism introduced in section~\ref{sec:inverse_uq}, the quantities of interest are represented by the 
map $\mathcal{G}(\theta)$, while the parameter vector $\theta \in \mathcal X$ defines a particular aortic shape.
The map is built on two recent works concerning diffeomorphic registration of aortic shapes~\cite{romor2026shape} and generative modeling of shapes based on stochastic interpolants~\cite{KATZ2026119120}.
In this section, we first introduce the geometry dataset used for training the registration algorithm and the generative model (section~\ref{ssec:geometry}), as well as for defining a latent parameter space encoding the main geometrical features (section~\ref{ssec:param-space}). In section~\ref{ssec:map} we detail the different steps of the forward map, entailing a generative model for the construction of a surface mesh, a smoothing step for the registration of a reference hexahedral mesh, the numerical solution of the Navier--Stokes equations, and the evaluation of an observation functional (e.g., a pressure drop between selected locations). The WFR method of section~\ref{subsec:wfr} requires
a differentiable surrogate map for the evaluation of the likelihood score. The steps to define this surrogate and 
the differences with the full forward map are discussed in section~\ref{sec:forward_lddmm}.
A sketch of the different steps defining the
forward map is provided in figure~\ref{fig:scheme_forward}.  %

\subsection{Geometry dataset and shape registration}\label{ssec:geometry}

The building blocks of the forward map rely on a dataset of synthetic aortic shapes created
using statistical shape modeling (SSM) from an initial cohort of $228$ subjects ($106$ aortic coarctation patients, $37$ post-treatment patients, and $85$ healthy individuals). The procedure
resulted in a dataset of $1257$ aortic shapes, described via a surface triangular mesh, 
a set of centerline control points (same amount for all shapes), 
and for each centerline point, the radius of the corresponding inscribed sphere.
For the detailed methodology, we refer the reader, e.g., to~\cite{goubergrits2022ct, thamsen2020unsupervised,thamsen2021synthetic,versnjak2024deep, yevtushenko2021deep}.

\begin{figure}[!h]
    \centering
    \includegraphics[width=0.3\textwidth]{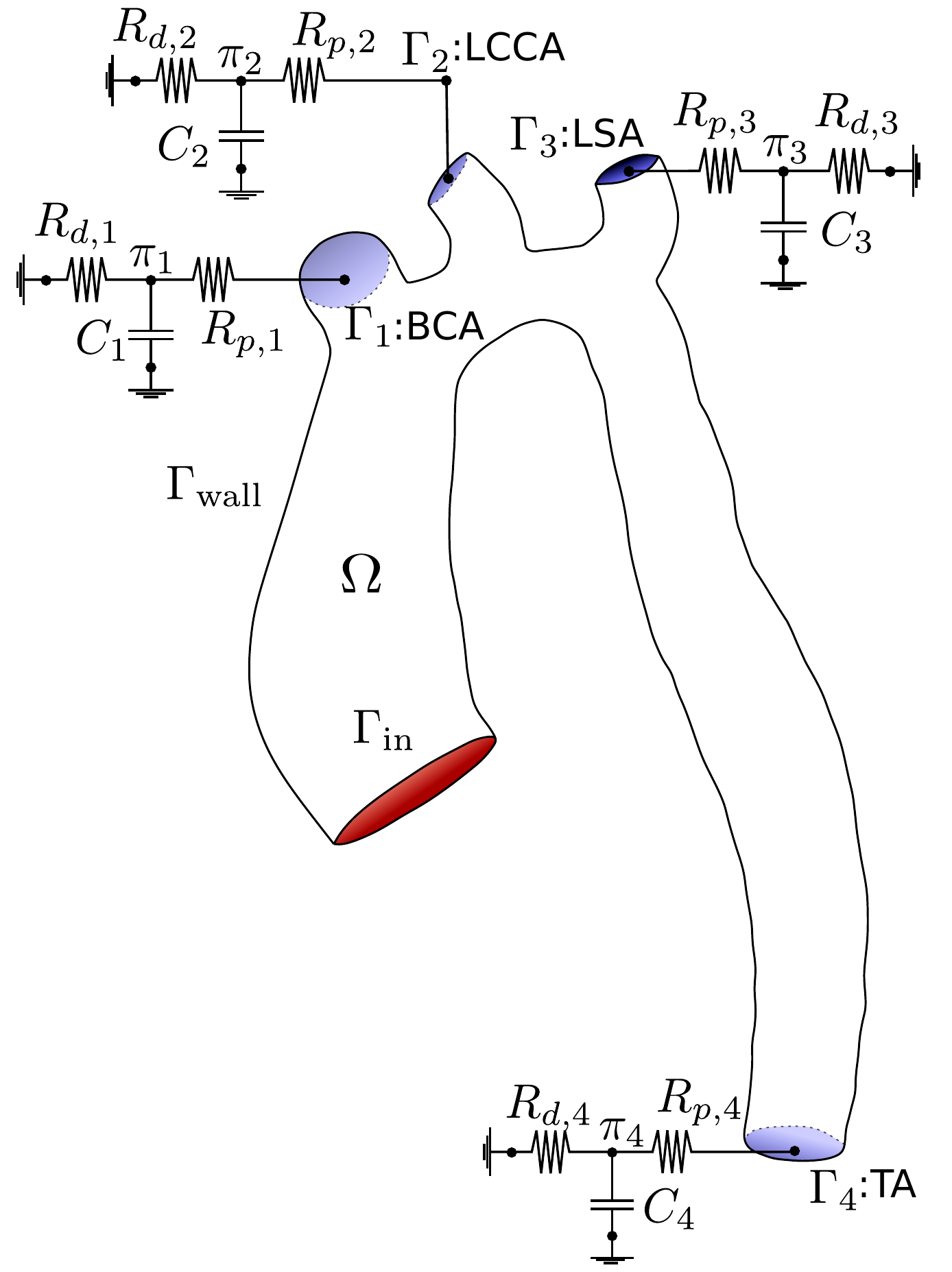}
    \includegraphics[width=0.69\textwidth]{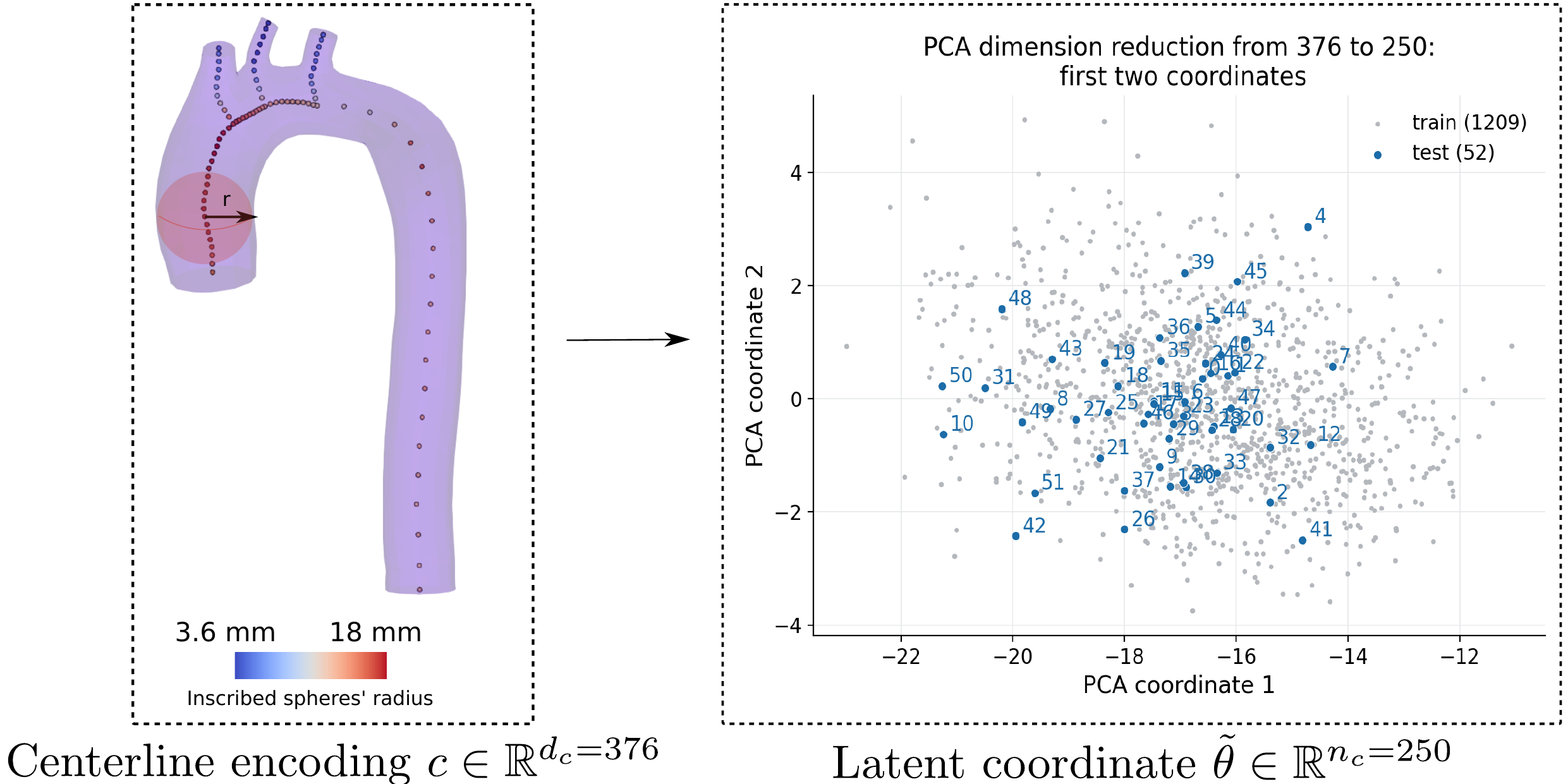}
    \caption{Left. Considered region of the aorta, highlighting the four outlet boundaries 
    (BCA: brachiocephalic artery, LCCA: left common carotid artery, LSA: left subclavian artery, TA: thoracic aorta).
    Right. Dimension reduction for the definition of the parameter space, encoding $94$ $4$-dimensional centerline 
    coordinates into a latent variable ($\tilde \theta$) of dimension $250$.}
    \label{fig:dataset}
\end{figure}

\begin{figure}[!ht]
  \includegraphics[width=0.9\textwidth]{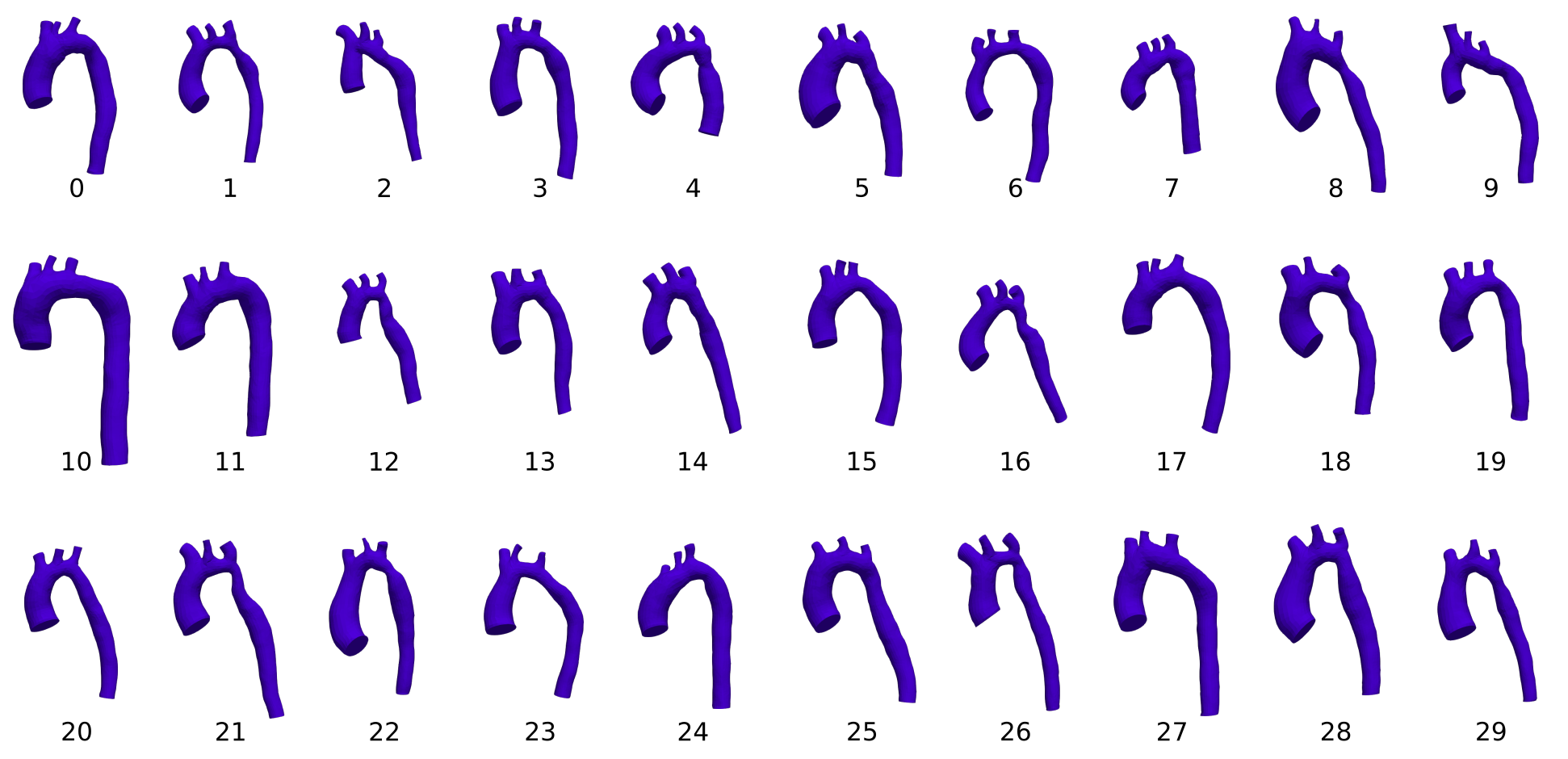}
  \caption{Selection of $30$ of the $52$ test geometries.}
  \label{fig:all_test_geometries}
\end{figure}

All shapes describe the same region of interest, comprising the vessel surface of the aortic arch up to the thoracic aorta (TA), including three main branches and the corresponding boundary surfaces (brachiocephalic artery, BCA, left common carotid artery, LCCA, left subclavian artery, LSA) (figure~\ref{fig:dataset}, left).
The dataset of $1257$ geometries is then split into $1205$ %
training and $52$ test shapes, as in~\cite{KATZ2026119120}. One training shape is selected as \textit{template} and indexed
with $i=0$. 
A sample of the test dataset is shown in figure~\ref{fig:all_test_geometries}.

Based on the SSM construction, each synthetic geometry is parametrized via a \textit{centerline encoding}
$c \in \mathbb{R}^{d_c}$, with $d_c:= 94 \cdot 4 = 376$, 
gathering the three coordinates and the radii of the corresponding inscribed spheres of $94$ selected control points along the centerline 
(main vessel and branches, figure~\ref{fig:dataset}, right).
For each shape in the database, we compute a bi-Lipschitz homeomorphic registration map on the template surface using the multi-level ResNet-LDDMM algorithm
proposed in~\cite{romor2026shape}. The template mesh is registered onto the point cloud describing other shapes, allowing us to define push-forward and pull-back
maps between the different computational domains and finite element solutions. This step is crucial to define the generative model introduced in section~\ref{sssec:lddmm}.

\begin{remark}
The formulation of the method does not depend on the particular choice of the template. 
However, since the template shape defines the reference space onto which the other shapes and numerical solutions are mapped, one should avoid
selecting shapes with extremely marked features.
For a more detailed discussion on the choice of the template we refer to~\cite{romor2026shape}.
\end{remark}

\subsection{Parameter space}\label{ssec:param-space}
\begin{figure}[htp!]
    \includegraphics[width=0.7\textwidth]{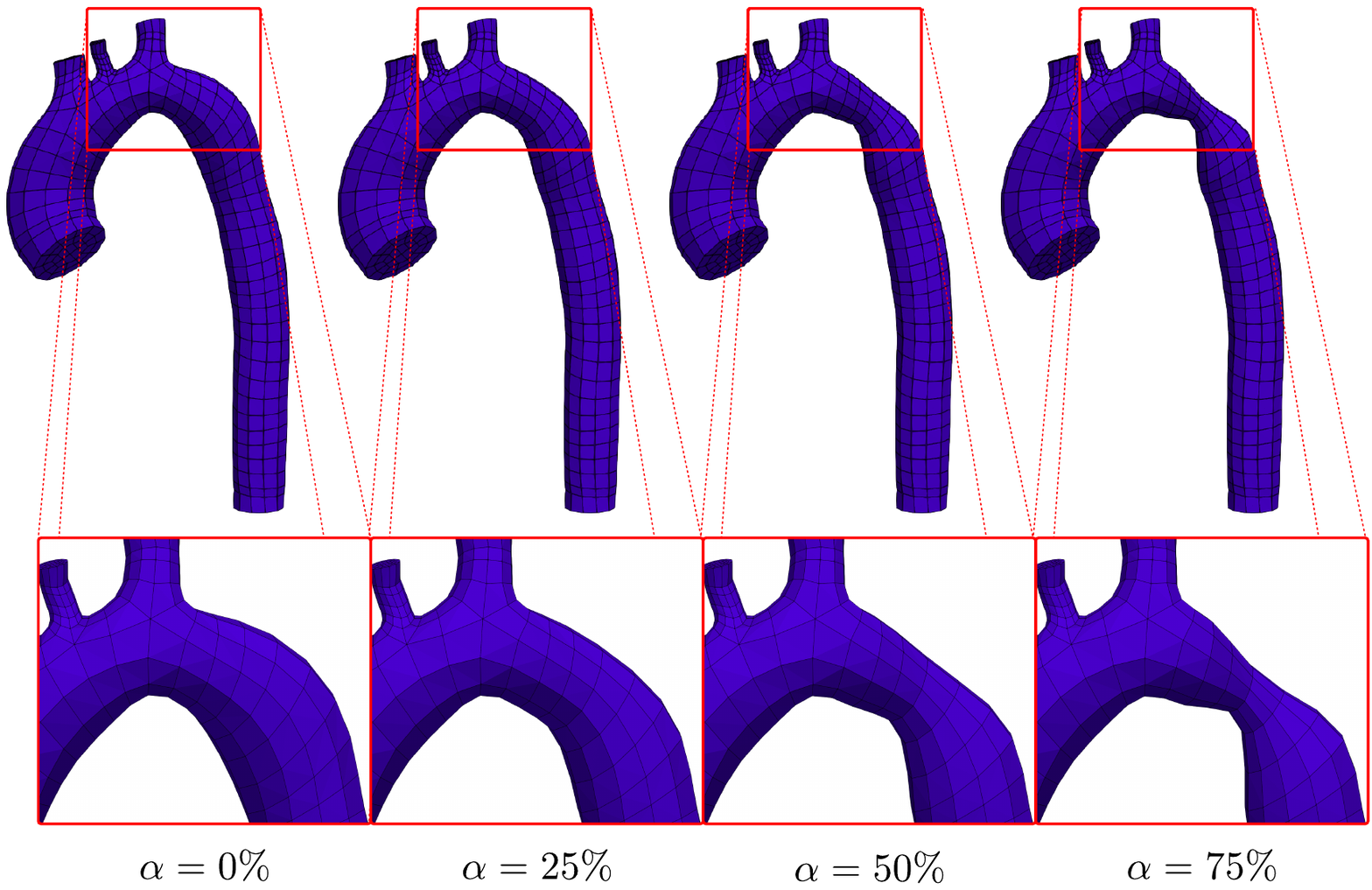}
    \caption{Parametric stenosis deformation for increasing severity levels $\alpha\in\{0\%, 25\%, 50\%, 75\%\}$ (from left to right).} %
    \label{fig:stenosis}
\end{figure}

The parameter space $\mathcal{X}$ is defined in two steps. 
First, we reduce the dimension of the centerline encoding using Principal Component Analysis (PCA) on the training set and retaining the first $n_c = 250$ modes.
Next, to further increase the variability of the considered population, we include an additional parameter describing the degree of a stenosis.
The stenosis is included via a one-dimensional parametric deformation, gradually reducing the radius of the inscribed sphere around a fixed location along the aorta centerline (see figure~\ref{fig:stenosis}).

The resulting parameter space $\mathcal X=\mathbb{R}^{n_\theta}$ has dimension $n_\theta=n_c+1$. For a parameter
$\theta = (\tilde{\theta},\tilde\alpha) \in \mathcal X$, the component $\tilde \theta \in \mathbb R^{n_c}$ will be referred to as the \textit{latent coordinate}, whilst
$\tilde \alpha \in [0,1]$ refers to the \textit{stenosis degree}.

Formally, let us denote by $\{c_i \in \mathbb R^{d_c}\}_{i=1}^{n_{\rm train}}$  the centerline encoding of the
training set, and with $U \in \mathbb{R}^{d_c\times n_c}$ the resulting PCA basis. Moreover, we 
denote by $m_c\in\mathbb{R}^{n_c}$ and $L_c\in\mathbb{R}^{n_c\times n_c}$ the mean and the square root of the covariance of the reduced representations of the training distribution $\{U^Tc_i\}_{i=1}^{n_{\text{train}}}\subset\mathbb{R}^{n_c}$, respectively. 
The prior distribution for the
parameter $\theta = (\tilde{\theta},\tilde\alpha)$  is then defined as a standard normal distribution on the space $\mathcal X$, after centering and scaling the latent coordinate with $m_c\in\mathbb{R}^{n_c}$ and $L_c\in\mathbb{R}^{n_c\times n_c}$.

\begin{remark}
The UQ methods introduced in section~\ref{sec:inverse_uq} require sampling
the parameter space. The dimension reduction, from $376$ to $250$ coordinates, has been implemented
with the purpose of retaining the most relevant geometrical features and limiting the dimension of the parameter space.
To this purpose, the truncation at $250$ was chosen by computing the PCA reconstruction error for a varying number of modes $50,100,\hdots,350$ and selecting the value at which the decay of the singular values shows the largest spectral gap. 
\end{remark}

\begin{remark}
    The location and the form of the additional stenosis are arbitrary. The purpose of including the stenosis in the parameter space is to demonstrate that any additional parametric deformation can be combined with the LDDMM stochastic interpolant used as generative prior for the surface.
\end{remark}

\subsection{Parameter-to-observable map}\label{ssec:map}

\subsubsection{Generative model for surface mesh}\label{sssec:lddmm}
The first step concerns the generation of a surface geometry. The approach
is based on the LDDMM stochastic interpolant framework introduced in~\cite{KATZ2026119120}
in the context of aortic blood flow simulations.

Given $\theta = (\tilde{\theta},\tilde \alpha)$, we first recover the corresponding centerline encoding $c(\theta) \in \mathbb R^{d_c}$ via 
\begin{equation}
\label{eq:latent_reconstruction}
        c(\theta)=U\big(L_c\tilde\theta+m_c\big).
\end{equation}
The vector $c(\theta)$ is used as a conditioning variable to define the LDDMM stochastic interpolants~\cite{KATZ2026119120}, which are obtained by solving an SDE of the form
    \begin{equation}\label{eq:stoc-interp}
    \begin{aligned}
        \mathrm{d}I_t & = b_t^{w}(I_t,c(\theta))\,\mathrm{d}t + \sigma_t\,\mathrm{d}W_t,\\
       \qquad I_{t=0}& =X_0.
       \end{aligned}
    \end{equation}
    In \eqref{eq:stoc-interp}, $X_0$ denotes
    the surface of the template shape and $b_t^{w}(I_t,c(\theta))$ is a drift coefficient
    trained using surface registration maps between the 
    shapes of the training dataset and the template~\cite{KATZ2026119120}. Moreover, $\sigma_t = \sqrt{t(1-t)}$ and $W_t$ is a standard Brownian motion.
    Equation \eqref{eq:stoc-interp} is solved using the forward Euler method with $10$ time steps up to $t=1$ and the resulting surface displacement is given by $\Psi_s = \mathbb E[I_{t=1}]$.

Next, we apply the deformation related to the stenosis parameter. 
The stenosis degree is obtained from the variable $\tilde \alpha \in \mathbb{R}$ 
through the standard-normal CDF, rescaled to a physical range, i.e., 
$$
\alpha=\alpha_{\min}+(\alpha_{\max}-\alpha_{\min})\,\cdfN(\tilde\alpha)\,.
$$
Based on the value of $\alpha$, the radii of the inscribed spheres are reduced gradually, 
up to a factor $\alpha$ at the selected location (see figure~\ref{fig:stenosis}), directly on the generated surface mesh.

\subsubsection{Mesh}\label{sssec:mesh}
The next step handles the extension of the surface deformation to a suitable hexahedral mesh of the bulk domain.
In particular, since the mesh is required for the numerical solution of the Navier--Stokes equations, it is necessary to guarantee an acceptable mesh regularity and avoid element degradation, especially in the highly deformed regions. The mesh generation process consists of two main steps.

First, the surface is preprocessed, applying a Taubin smoothing with the graph Laplacian and correcting the inlet and outlet faces by projecting the mesh nodes on planar surfaces.

The hexahedral mesh in the bulk domain is defined by mapping the template volume mesh onto the new surface.
We solve a multi-step linear elasticity problem for an internal displacement field, starting from the template mesh and gradually deforming it towards the generated surface $\Psi_s$:
\begin{equation}
  \label{eq:mesh_extension}
    \left\{
    \begin{aligned}
        \nabla\cdot\!\left[\tfrac{E(\mathbf{x})}{2(1+\nu_{\mathrm{P}})}\big(\nabla\Psi_n+\nabla\Psi_n^{T}\big)+\tfrac{E(\mathbf{x})\,\nu_{\mathrm{P}}}{(1+\nu_{\mathrm{P}})(1-2\nu_{\mathrm{P}})}(\nabla\cdot\Psi_n)\,\mathrm{Id}\right] &= 0 &&\text{in }\Psi_{n-1}(\Omega_0),\\
        \Psi_{n} &= \tfrac{e^{-n/N_{\mathrm{steps}}}-1}{e^{-1}-1}\,\Psi_s &&\text{on }\partial\Psi_{n-1}(\Omega_0).
    \end{aligned}
    \right.
\end{equation}
In \eqref{eq:mesh_extension}, $\Omega_0$ denotes the template domain, 
$N_{\mathrm{steps}}=8$ is the number of pseudo-time steps, and $\Psi_s$ is the surface generated from \eqref{eq:stoc-interp}, including, for ease of notation, also the smoothing and planar correction steps mentioned above.
The Poisson ratio is defined as $\nu_{\mathrm{P}}=0.3$, and the stiffness $E(\mathbf{x})$ is set cell-wise to the largest ratio between the cell diameter and the length of one of its edges, so that already distorted cells are stiffer and deform less than well-shaped ones. 
The procedure is designed to gradually deform the template volume mesh onto a mesh compatible with the new surface mesh, preserving the original mesh topology. 
However, it does not guarantee that elements do not degenerate after each iteration. 
Therefore, after each pseudo-time step $n$, the cells with non-positive Jacobian determinant or with aspect ratio above $100$ are repaired by at most $40$ gradient-descent steps on a cell aspect-ratio functional. We refer to~\cite{KATZ2026119120} for the details of the mesh generation procedure.
In what follows, we will denote by $\Omega_{\theta}$ the computational domain generated for the parameter $\theta \in \mathcal X$.
    
\begin{remark}[Handling of invalid geometries]
\label{rmk:invalid}
The forward map $\mathcal{G}$ is defined on the whole parameter space $\mathbb{R}^{n_\theta}$. However,  
not every parameter vector corresponds to a valid geometry, since the prior distribution is a standard normal on the PCA space of the geometric conditioning vector $c \in \mathbb R^{d_c}$ (section~\ref{ssec:param-space}), which is only an approximation of the true distribution of admissible geometries. 
As a result, the map can produce invalid geometries, i.e., shapes with self-intersecting surfaces, with negative inscribed radii, or surfaces that cannot be meshed with the multi-step linear elasticity extension \eqref{eq:mesh_extension}.
In our study, when the multi-step linear elasticity extension fails to produce a valid mesh at a step $n<N_{\rm steps}$, we stop the iteration at that stage.
If the extension fails at the first step, the forward map returns the template shape. Invalid parameters are thus mapped to an interpolation between the template and the generated surface, so that the forward map always returns an admissible domain and no prior mass is mapped to invalid geometries.
\end{remark}

Notice that the numerical solver (section~\ref{sssec:cfd}) requires a hierarchy of meshes for the geometric multigrid preconditioner. Hence, the mesh extension \eqref{eq:mesh_extension} is applied recursively to the entire hierarchy generated from the template mesh.

\subsubsection{Numerical solution}\label{sssec:cfd}
On the computational domain $\Omega_{\theta}$, we model the blood as an incompressible Newtonian fluid, described by a velocity field $\bu: \Omega_{\theta} \times [0,T] \to \mathbb{R}^3$ and a pressure field $p: \Omega_{\theta} \times [0,T] \to \mathbb{R}$, which satisfy the Navier--Stokes equations
 \begin{equation}\label{eq:nse}
        \left\{
        \begin{aligned}
            \rho_{\mathrm{f}}\,\partial_t \bu + \rho_{\mathrm{f}}\,(\bu\cdot\nabla)\bu -
            2 \nabla \cdot \left( \mu_{\mathrm{f}}\mathbb D(\bu)\right) 
            + \nabla p &= 0, && \text{in }\Omega_{\theta},\\[2pt]
            \nabla\cdot \bu &= 0, && \text{in }\Omega_{\theta}\,.
        \end{aligned}
        \right.
    \end{equation}
    In \eqref{eq:nse}, $\mathbb D(\bu)$ is the symmetric part of the deformation tensor, 
    $\rho_{\mathrm{f}}$ is the blood density, and $\mu_{\mathrm{f}}$ the dynamic viscosity, which includes
    a turbulent viscosity contribution based on a $\sigma$-model~\cite{Nicoud2011}. 
    We consider a cardiac period $T=1\,\mathrm{s}$. Notice that, for ease of notation, we omitted the dependence on $\theta$ of the solution $(\bu,p)$. 

The system \eqref{eq:nse} is completed by the boundary conditions
 \begin{equation}\label{eq:nse-bc}
\left\{
\begin{aligned}
    \bu &= \bu_{\rm in} && \text{on }\Gamma_{\theta,\mathrm{in}},\\
    \bu &= 0 && \text{on }\Gamma_{\theta,\mathrm{wall}},\\
    \big(-p\,\mathrm{Id}+\mu_{\mathrm{f}}(\nabla \bu+\nabla \bu^{T})\big)n &= -P_i\,n && \text{on }\Gamma_{\theta,i},\ \ i\in\{1,2,3,4\}.
\end{aligned}
\right.
\end{equation}
At the inlet $\Gamma_{\theta,\rm in}$, we impose a given profile in the form $\bu_{\rm in}(t,\bx) = a(t)\,\beta(\bx)$, obtained from a reference spatial profile 
$\beta(\bx)$, registered to the current inlet shape, and multiplied by a time-periodic function $a(t)$ \cite{romor2026shape}. On the vessel wall 
$\Gamma_{\theta,\mathrm{wall}}$ we impose
no-slip Dirichlet boundary conditions, i.e., neglecting the interaction with the vessel wall.
On the four outlet surfaces $\Gamma_{\theta,i}$ we impose
an average pressure $P_i$ obtained via a so-called 
three-element Windkessel model, i.e., as the solution of the ODE
\begin{equation}\label{eq:wk}
C_{d,i}\frac{\mathrm{d}}{\mathrm{d}t}(P_i-R_{p,i}\,Q_i)+\frac{P_i-R_{p,i}\,Q_i}{R_{d,i}}=Q_i.\end{equation}
Equation \eqref{eq:wk} depends on the outgoing flow $Q_i$ on the outlet $\Gamma_{\theta,i}$ and on three parameters, $R_{p,i}$, $R_{d,i}$, and
$C_{d,i}$, which model the resistance and capacitance of the vascular system downstream.
The Windkessel parameters are calibrated
across the different shapes as described
in~\cite{romor2026shape}, depending on 
a given flow split among the four outlets.

The system of equations \eqref{eq:nse}-\eqref{eq:nse-bc} is discretized using a
dual-splitting scheme in time and a
high-order discontinuous Galerkin (DG2-DG1) discretization in space, with matrix-free geometric multigrid preconditioners in \texttt{ExaDG}~\cite{ExaDG2020, KATZ2026119120}
with AVX SIMD vectorization and MPI parallelism on $32$ cores. 
The geometric multigrid preconditioner requires a hierarchy of meshes. The finest mesh is composed of $2271$ nodes and $1792$ hexahedra, resulting
in $159488$ spatial degrees of freedom.

\subsubsection{Observations}\label{sssec:observations}
The last step in evaluating? the forward map $\mathcal G(\theta)$ is defined by evaluating an observation operator $\mathcal{O}(\bu,p)$ on the numerical solution of \eqref{eq:nse}-\eqref{eq:nse-bc}. 

For the examples considered in this work, we define a subset of $42$ centerline points (shown in figure~\ref{fig:scheme_forward} (observation operator) on the aortic arch of a sampled geometry) 
with corresponding cross-sections $\{S_i\}_{i=1}^{42}$, and 
consider the following observations.
\begin{enumerate}
    \item The pressure drop along the centerline at $t=0.25\,\mathrm{s}$ with respect to the average inlet pressure, defined as
    \begin{equation}
      \label{eq:obs_1}
        \mathcal{O}_{1}(\bu,p)
        =\mathcal{O}_{1}(p) = \left(\frac{1}{|\Gamma_{\theta,\rm in}|}\int_{\Gamma_{\theta,\rm in}}
\!\!\left(p\left(t=0.25\,\mathrm{s},\bx\right)-\frac{1}{|S_i|}\int_{S_i}
        p\left(t=0.25\,\mathrm{s},\bx\right)\right)\,\mathrm{d}x\right)_{i=1}^{42}.
    \end{equation}
    \item The maximum pressure drop across all $42$ centerline sections:
    \begin{equation}
      \label{eq:obs_2}
        \mathcal{O}_{2}(\bu,p)=\mathcal{O}_{2}(p)=\max_{i=1,\dots,42} \left|\mathcal{O}_{1}(p)\right|_i.
    \end{equation}
\end{enumerate}
Notice that we select the time instant $t=0.25\,\mathrm{s}$ as it lies right after the systolic peak.

We focus on observations only dependent on the pressure field. However, the framework can be analogously applied to other quantities of interest (centerline velocities, secondary-flow indices, wall shear stresses, etc.).

The steps of the forward map described in the previous sections are summarized
in algorithm~\ref{alg:g-theta} and figure~\ref{fig:scheme_forward}.

\begin{algorithm}[t]
\caption{Forward map} %
\label{alg:g-theta}
\begin{algorithmic}[1]
\itemsep=2pt
\Function{ForwardMap}{$\theta = (\tilde \theta,\tilde \alpha)$} \Comment{latent coordinate and stenosis degree}
\State latent reconstruction $c = U\big(L_c\tilde\theta + m_c\big)$ %
\State generate surface mesh via LDDMM stochastic interpolants \eqref{eq:stoc-interp}
\State include stenosis defined by $\tilde \alpha$
\State generate volume mesh (smoothing and multi-step elasticity \eqref{eq:mesh_extension})
\State compute $(\bu,p)$ solution to \eqref{eq:nse}-\eqref{eq:nse-bc}
\State evaluate observations $\mathcal G(\theta) = \mathcal O(p)$
\State \Return $\mathcal G(\theta)$
\EndFunction
\end{algorithmic}
\end{algorithm}

\subsection{Geometry-conditioned autoregressive LDDMM stochastic interpolant}
\label{sec:forward_lddmm}
The WFR method (section~\ref{subsec:wfr}) requires a differentiable surrogate forward map $\Tilde{\mathcal G}$  to evaluate the likelihood score $-\nabla_\theta\slike(\theta)$ in~\eqref{eq:wfr_langevin}. 

To this end, we replace the high-fidelity PDE solver (section~\ref{sssec:cfd}) with a machine learning surrogate that emulates the pressure field over the cardiac cycle, conditioned on the geometry given by the centerline encoding vector.

The pressure surrogate is a \emph{flow-matching} (stochastic-interpolant) model~\cite{KOHL2026108641} of the transient lumen pressure field of a patient-specific aorta, and it relies heavily on the multi-step LDDMM registration of~\cite{romor2026shape} to produce the training dataset of pressure fields supported on the same template. 

The surrogate is trained to predict the pressure field $p$ on the whole domain. To increase the resolution of the predictions, the training pressure fields are interpolated onto a DG2 finite element space on the same mesh.
Hence, conditioned on the geometry, the surrogate produces the scalar pressure $p$ at each of the $n_{\text{nodes}}=48\,384$ support points of the DG2 space
for the physical time instants $t\in\{0,0.05,0.1,0.15,0.2,0.25\}\,\mathrm{s}$ within the cardiac period. The SDE describing the model reads
\begin{equation}
  \label{eq_lddmm_p}
  \begin{aligned}
    \mathrm{d}I^p_t & = b^p_t(I^p_t,\mathbf{x}, c)\,\mathrm{d}t + \sigma^p_t\,\mathrm{d}W_t,\\
    I^p_{t=0}&=p_0,
    \end{aligned}
\end{equation}
where $\bx \in\mathbb{R}^{n_{\text{nodes}}\times 3}$ are the coordinates of the mesh vertices
and $c$ the centerline encoding.

Rather than predicting each snapshot directly, the network learns the time-dependent drift $b^p$ that transports one snapshot into the next. For a source/target pair of snapshots at times $t$ and $t+\Delta$, the interpolant is defined as
$$
I^p(s)=(1-s)\,p_t+s\,p_{t+\Delta}+\sigma^p(s)\,\xi,
$$
at interpolation time $s\in[0,1]$, with Karras scheduler $\sigma^p(s)=\sigma_{\text{max}}\left(\sigma_{\text{min}}/\sigma_{\text{max}}\right)^{s^{\rho}}$~\cite{karras2022elucidating} and $\xi\sim\mathcal N(0,\mathrm{Id})$. The neural network regresses the target drift $(p_{t+\Delta}-p_t)+\dot\sigma^p(s)\,\xi$. 
At inference, the learned drift is integrated by the explicit Euler method with $10$ substeps per interval to roll the pressure field forward across the instants $t\in\{0,0.05,0.1,0.15,0.2,0.25\}\,\mathrm{s}$.

At each evaluation, the drift network receives, for each node, 
\begin{itemize}
    \item the current pressure state (scalar),
    \item the node coordinates $\mathbf{x} \in\mathbb{R}^{n_{\text{nodes}}\times 3}$, lifted to Fourier features, and
    \item a single \emph{global conditioning token} built from the geometry conditioning vector $c(\theta) \in\mathbb{R}^{d_c}$ (the centerline control points and inscribed radii) together with two sinusoidal time embeddings (global cycle time and local interpolation time),
\end{itemize} 
and outputs one scalar per node, i.e., the pressure drift field.

The node features are reordered along a Morton (Z-order) space-filling curve and packed into a dense $3$D voxel grid. The ordering does not reproduce the connectivity of the mesh, but it makes the fields available on a regular grid, on which convolutional and transformer blocks are computationally more efficient than on the unstructured mesh, while keeping \emph{spatially} close nodes close in memory. %
A depth-$4$ U-Net of \emph{conditioned transformer blocks} processes the volume, with channel widths $64\!\to\!128\!\to\!256\!\to\!512$ and a $1024$-channel bottleneck. The output is un-permuted back to the nodes and projected to the scalar drift. The reported architecture uses $6$ transformer blocks per stage and $6$ bottleneck blocks ($\approx\!112$M parameters, see table~\ref{tab:model-params}).

\begin{table}[t]
  \centering
  \caption{Parameter count of the stochastic interpolant drift network.}
  \label{tab:model-params}
  \definecolor{hdrblue}{RGB}{31,78,121}
  \definecolor{rowblue}{RGB}{222,235,247}
  \definecolor{totalblue}{RGB}{189,215,238}
  \begin{tabular}{lrr}
    \toprule
    \rowcolor{hdrblue}
    {\color{white}\textbf{Component}} & {\color{white}\textbf{Parameters}} & {\color{white}\textbf{Share}} \\
    \midrule
    Fourier PE (position $+$ pressure)              & 2,128       & $<0.1\%$ \\
    Conditioning embedding (geometry $+$ 2$\times$time) & 98,624   & $0.1\%$  \\
    U-Net input embed $+$ conditioning proj.        & 106,880     & $0.1\%$  \\
    Encoder 4 stages (ch.\ $64\!\to\!512$)      & 19,453,888  & $17.4\%$ \\
    Bottleneck (ch.\ $1024$; conv $+$ attention)    & 63,507,456  & $56.8\%$ \\
    Decoder 4 stages (ch.\ $512\!\to\!64$)      & 28,683,200  & $25.6\%$ \\
    Output (norm $+$ linear $\to 1$)                & 193         & $<0.1\%$ \\
    \midrule
    \rowcolor{totalblue}
    \textbf{Total} & \textbf{111,852,369} & \textbf{100\%} \\
    \bottomrule
  \end{tabular}
\end{table}

For the surrogate forward map, we restrict to the case
$\tilde \alpha = 0$ (without additional stenosis) degree, considering only the parameters
of the latent space describing the geometric encoding (centerline and inscribed radii), i.e., $\theta=(\tilde\theta,0)$. 
This restriction is necessary, since the surrogate
model is trained on numerical solutions obtained with the original training dataset, which contains the shapes without the stenosis. An extension with 
an additional parametric deformation would require the generation of a new training dataset.
The strong dependence on the training distribution is one important limitation of this surrogate approach, 
and retraining or fine-tuning is needed whenever the training distribution must be updated.
A further difference with respect to the forward map introduced in section~\ref{ssec:map}
is that the volume mesh is generated using radial basis function (RBF) interpolation 
to extend the registered surface displacement into the interior of the mesh, instead of using the multi-step elasticity \eqref{eq:mesh_extension}. 
RBF interpolation is needed to make the forward map amenable to batched execution on GPUs. 
Notice that the surrogate map does not require smoothing steps to improve the mesh quality, since the numerical solution of the Navier--Stokes equation, which
might be sensitive to the computational mesh, is no longer needed.

\begin{figure}[htp!]
    \includegraphics[width=.8\textwidth]{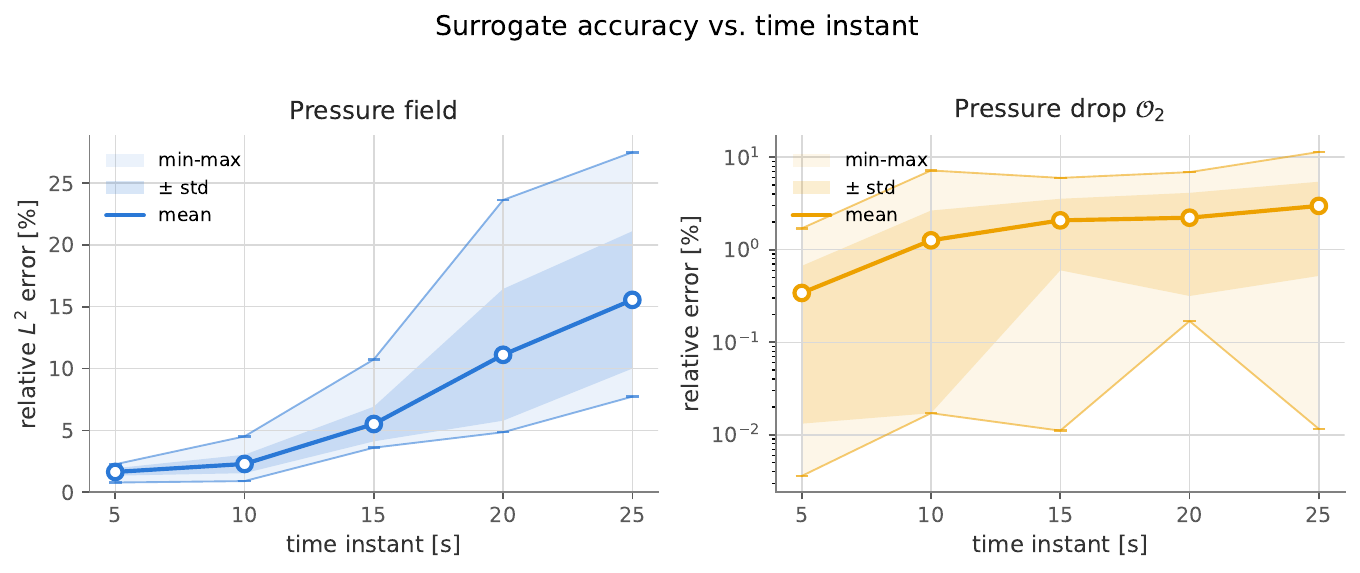}
    \caption{LDDMM stochastic interpolant pressure surrogate: $l^2$-relative error (left) and relative error on the maximum pressure drop (right), i.e., the output of the observation operator~\eqref{eq:obs_2}.}
    \label{fig:surrogate_error}
\end{figure}

\begin{figure}[htp!]
\centering
    \includegraphics[width=.8\textwidth]{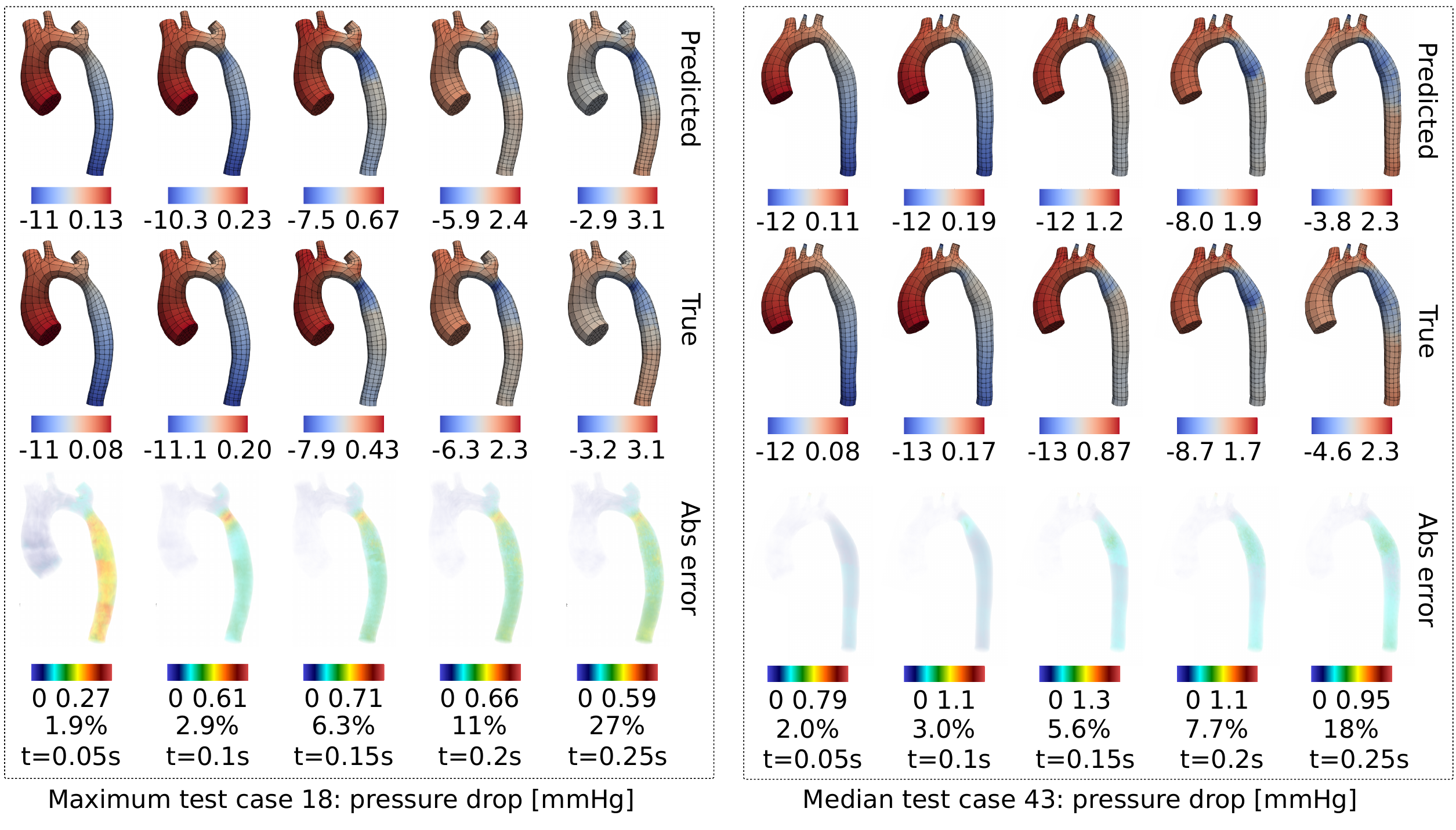}
    \caption{LDDMM stochastic interpolant pressure surrogate: snapshots of the cases with the maximum $l^2$-relative error (left, case $18$) and the median $l^2$-relative error (right, case $43$), averaged over all time instants.}
    \label{fig:surrogate}
\end{figure}

Figure~\ref{fig:surrogate_error} reports the relative $l^2$ error on the pressure field and the relative error on the maximum pressure drop for the $52$ test cases. Predicted and high-fidelity pressure snapshots are shown in figure~\ref{fig:surrogate} for the cases with the maximum and the median $l^2$-relative error.

\begin{figure}[htp!]
    \includegraphics[width=0.7\textwidth]{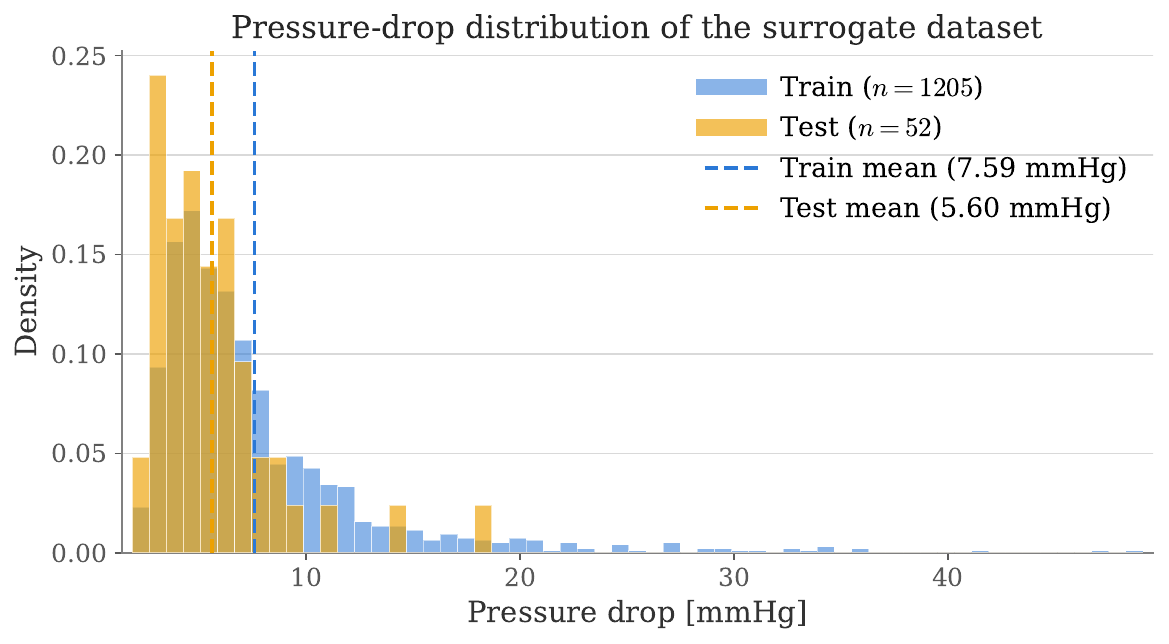}
    \caption{Pressure drop distribution of the training and test datasets.}
    \label{fig:surrogate_dist}
\end{figure}

The distribution of the pressure drop obtained from the surrogate 
for the training and test datasets (without stenosis) is shown in figure~\ref{fig:surrogate_dist}. 

\begin{figure}[htp!]
  \includegraphics[width=\textwidth]{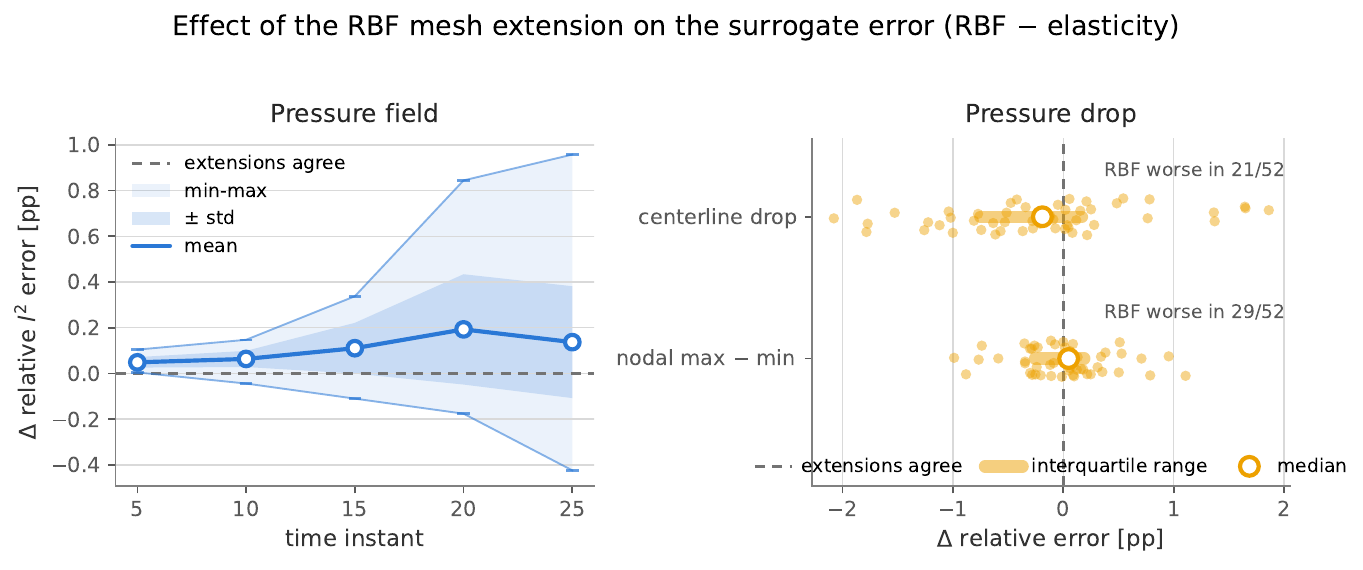}
  \caption{Effect of the RBF interpolation on the accuracy of the surrogate, measured against the same high-fidelity solution. Left: difference between the $l^2$ relative errors on the pressure field, at each time instant, with mean, $\pm$ standard deviation band and min-max range over the $52$ test cases. Right: difference between the relative errors on the centerline pressure drop $\mathcal{O}_1$ and on the maximum pressure drop $\mathcal{O}_2$. Each marker corresponds to one test case, with median and interquartile range. Both left and right charts are in percentage points (pp): on the left, values of $1$ correspond to $1\%$ variation of the relative $l^2$ error.}
  \label{fig:rbf_sensitivity}
\end{figure}

Notice that, as a result of the RBF interpolation, the surrogate is evaluated on volumetric meshes that differ from the ones it has been trained on. 
The errors in figures~\ref{fig:surrogate_error} and \ref{fig:surrogate} are measured on the meshes produced by the multi-step elasticity extension. 
To quantify the effect of this substitution on the accuracy of the surrogate, we compare the two options case by case against the same high-fidelity solution.
Figure~\ref{fig:rbf_sensitivity} shows that the RBF interpolation affects the relative errors on the pressure drop and on the observed quantities $\mathcal{O}_1$ and $\mathcal{O}_2$ by at most $1\%$. Thus the surrogate is only mildly sensitive to the choice of the mesh extension approach.
The surrogate forward map $\Tilde{\mathcal{G}}$ used by the WFR method
is summarized in algorithm~\ref{alg:g-tilde-theta} and figure~\ref{fig:scheme_forward}.

\begin{algorithm}[ht]
\caption{Surrogate forward map} %
\label{alg:g-tilde-theta}
\begin{algorithmic}[1]
\itemsep=2pt
\Function{SurrogateForwardMap}{$\theta = (\tilde \theta,\tilde \alpha=0)$} \Comment{latent coordinate only}
\State latent reconstruction $c = U\big(L_c\tilde\theta + m_c\big)$ %
\State generate surface mesh via LDDMM stochastic interpolants \eqref{eq:stoc-interp}
\State RBF interpolation to extend the surface mesh to the volume 
\State $p$ = flow-matching pressure surrogate \eqref{eq_lddmm_p}
\State evaluate observations $\Tilde{\mathcal G}(\Tilde\theta) = \mathcal{O}(p)$ 
\State batched evaluation of $\nabla_{\Tilde\theta}\Tilde{\mathcal G}(\Tilde\theta)$ by backward automatic differentiation
\State \Return $\{\Tilde{\mathcal G}(\Tilde{\theta})\},\{\nabla\Tilde{\mathcal G}(\Tilde{\theta})\}$
\EndFunction
\end{algorithmic}
\end{algorithm}

\begin{remark}[Approximation of the gradient of the surrogate output with respect to the parameters]
\label{rmk:grad}
At each discrete time instant in which the pressure is inferred, the LDDMM stochastic interpolant takes as input the coordinates of the nodes of the generated mesh, $\mathbf{x}(\theta)$, and the centerline encoding $c(\theta)$, which is used as input to the conditioning branch. 
  To approximate the gradient of the surrogate $I^p_{t=0.25\,\mathrm{s}}(\mathbf{x}(\theta),c(\theta))$, we neglect the partial derivative with respect to the mesh nodes and consider a 
  recursive gradient evaluation:
  \begin{equation*}
    \frac{\partial I^p_{t_k}}{\partial \theta} = \frac{\partial I^p_{t_{k-1}}}{\partial \theta} + \Delta s \left(\frac{\partial b^{p}_{t_k}(I^p_{t_{k-1}}, \mathbf{x}, c)}{\partial I^p_{t_{k-1}}} \frac{\partial I^p_{t_{k-1}}}{\partial \theta} + \cancel{\frac{\partial b^{p}_{t_k}(I^p_{t_{k-1}}, \mathbf{x}, c)}{\partial \mathbf{x}}\frac{\partial\mathbf{x}}{\partial \theta}}+ \frac{\partial b^{p}_{t_k}(I^p_{t_{k-1}}, \mathbf{x}, c)}{\partial c}\frac{\partial c}{\partial \theta}\right),
  \end{equation*}
  so the approximate gradient reduces to:
  \begin{equation*}
    \frac{\partial}{\partial\theta} 
    I^p_{t=0.25\,\mathrm{s}}(\mathbf{x}(\theta),c(\theta)) \approx \Delta s\sum_{k=1}^{\lfloor 1/\Delta s\rfloor } \left[\prod_{m=k+1}^{\lfloor 1/\Delta s\rfloor}\left(\mathrm{Id}+\Delta s \frac{\partial b^{p}_{t_m}(I^p_{t_{m-1}}, \mathbf{x}, c)}{\partial I^p_{t_{m-1}}}\right)\right]\frac{\partial b^{p}_{t_k}(I^p_{t_{k-1}}, \mathbf{x}, c)}{\partial c}\frac{\partial c}{\partial \theta}.
  \end{equation*}
\end{remark}

\section{Numerical results}
\label{sec:results}
This section presents the numerical results
of the aortic shape inversion problem for the two selected types of observables (centerline and scalar pressure drops). 
We compare the chosen methods
in terms of convergence, posterior distributions of the parameters and observables, and representative aortic geometries from the posterior parameter distributions.

The EAKI (section~\ref{subsec:eaki}) is used
for the case of centerline pressure drop (operator defined in~\eqref{eq:obs_1}), as the observed/specified values of the pressure drop along the vessel centerline are informative enough to reconstruct a single geometry, while aBUS (section~\ref{subsec:abus}) and WFR (section~\ref{subsec:wfr}) are used for
the scalar pressure drop~\eqref{eq:obs_2}, for sampling multiple geometries associated with the same observed/prescribed value of the maximum pressure drop along the vessel centerline.

For EAKI and aBUS, the forward map is based on the high-fidelity PDE solver and on the mesh extended from the template (section~\ref{sec:forward}). 
For WFR, which additionally requires the gradient
$\frac{\partial \mathcal{G}}{\partial \theta}$, we use the surrogate map described in section~\ref{sec:forward_lddmm}.
The high-fidelity solver for the forward map is run in parallel for each particle of the ensembles on an HTCondor queue with $4$ nodes having $72, 128, 64, 128$ cores,
with each run requiring $32$ CPU cores. The surrogate model runs in batches of $16$ particles (forward and gradient evaluation) on a single NVIDIA H200 NVL4 141GB GPU.

In all cases, the prior distribution is the standard normal distribution in the latent parameter space, and the observation noise is relative, $\Sigma_{\eta}=\mathrm{diag}\big((0.1\,y)^2\big)$, i.e., $10\%$ of the assimilated value. 
The stenosis degree ($\alpha$) is not considered in the surrogate model, as described in section~\ref{sec:forward_lddmm}.
The setup of the three methods and the key differences are summarized in table~\ref{tab:methods_setup}.

\begin{table}[htp!]
  \centering
  \caption{Key properties of the three considered inversion methods.}
  \label{tab:methods_setup}
  \definecolor{hdrblue}{RGB}{31,78,121}
  \begin{tabular}{lccccc}
    \toprule
    \rowcolor{hdrblue}
    {\color{white}\textbf{Method}} & {\color{white}\textbf{Observation}}
      & {\color{white}\textbf{Stenosis}} & {\color{white}\textbf{Surrogate}}
      & {\color{white}\textbf{Gradient}} & {\color{white}\textbf{Gaussian posterior}}\\
    \midrule
    EAKI & centerline pressure drop~\eqref{eq:obs_1} & yes & no  & no  & yes  \\
    aBUS & pressure drop~\eqref{eq:obs_2} & yes & no  & no  & no  \\
    WFR  & pressure drop~\eqref{eq:obs_2} & no  & yes & yes & no \\
    \bottomrule
  \end{tabular}
\end{table}

\subsection{Centerline pressure drop}
\label{subsec:results_aortic_eaki}
We study the convergence behavior and the performance of the EAKI considering pressure data along the centerline.

To this end, we focus on 
two particular cases, i.e., the worst and the median test shapes based on the error of the pressure surrogate (cases $18$ and $43$, respectively, see figure~\ref{fig:surrogate}), with 
an additional stenosis with degree $\alpha=0.7$. With a little abuse of notation, we will refer to test case 18 and 43 in the  following part also to indicate the
geometries modified with the stenosis.
For these two cases, we solve the Navier--Stokes equations \eqref{eq:nse} and evaluate
the centerline pressure drop, defining the synthetic
measurements as $y=\mathcal{O}_1(p)$ at $t=0.25\,\mathrm{s}$ as defined in~\eqref{eq:obs_1}. The vector $y$ is then used as the observation for the EAKI method.

It is important to observe that the goal is not to reconstruct the original shapes (and, in general, this problem does not have a unique solution), but rather to obtain an ensemble of shapes -- in the EAKI case, approximating a Gaussian measure -- whose observation is close to the target one.

The most relevant parameters for the EAKI method are
\begin{itemize}
    \item the maximum number of iterations ($T_{\text{max}}$),
    \item the ensemble size ($N$),
    \item the time step ($\Delta\tau$), and
    \item the jitter threshold ($\sigma_{\text{thr}}$).
\end{itemize}
In the numerical tests, we vary the ensemble size and the jitter threshold. At the first iteration, the EAKI ensemble particles are sampled from Latin hypercube sampling inside $[-3, 3]^{n_\theta}$.
The number of iterations needed is usually small (below $10$), and we set $\Delta\tau=1/2$. For parameter studies related to $\Delta\tau$ and its interpretation, we refer to the original reference~\cite{huang2022efficient}.

\begin{figure}
  \includegraphics[width=\textwidth]{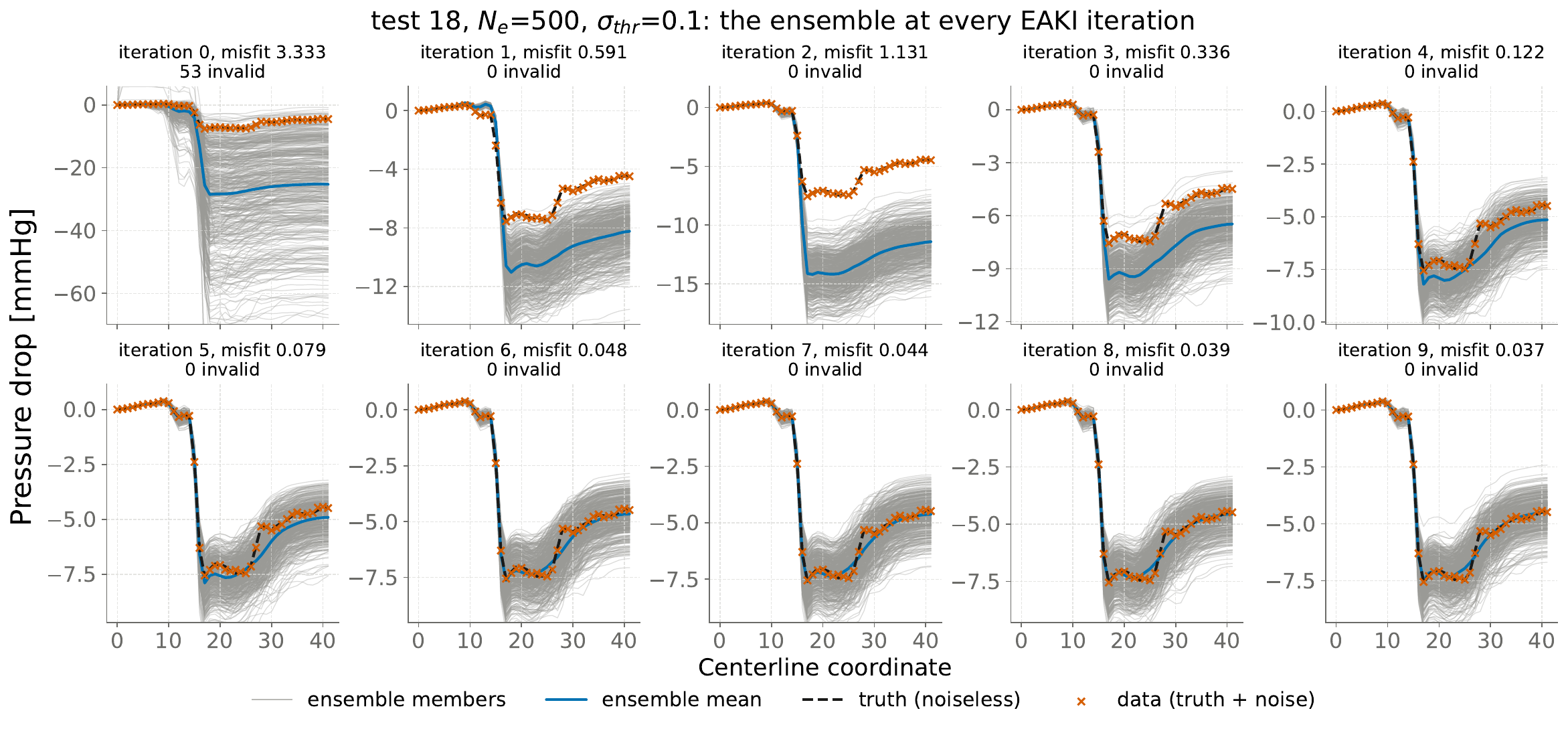}
  \caption{Test case $18$. Centerline pressure drops for the EAKI ensemble over $T_{\text{max}}=10$ iterations. Jittering threshold 
  $\sigma_{\text{thr}}=0.1$.}
  \label{fig:eaki18}
\end{figure}

We show the convergence of EAKI for test case $18$ in figure~\ref{fig:eaki18} for $\sigma_{\text{thr}}=0.1$.
The EAKI delivers, as a final result, an average shape given by the ensemble mean. The variance of the ensemble can be arbitrarily controlled through the jitter threshold $\sigma_{\text{thr}}$.

\begin{figure}
  \includegraphics[width=\textwidth]{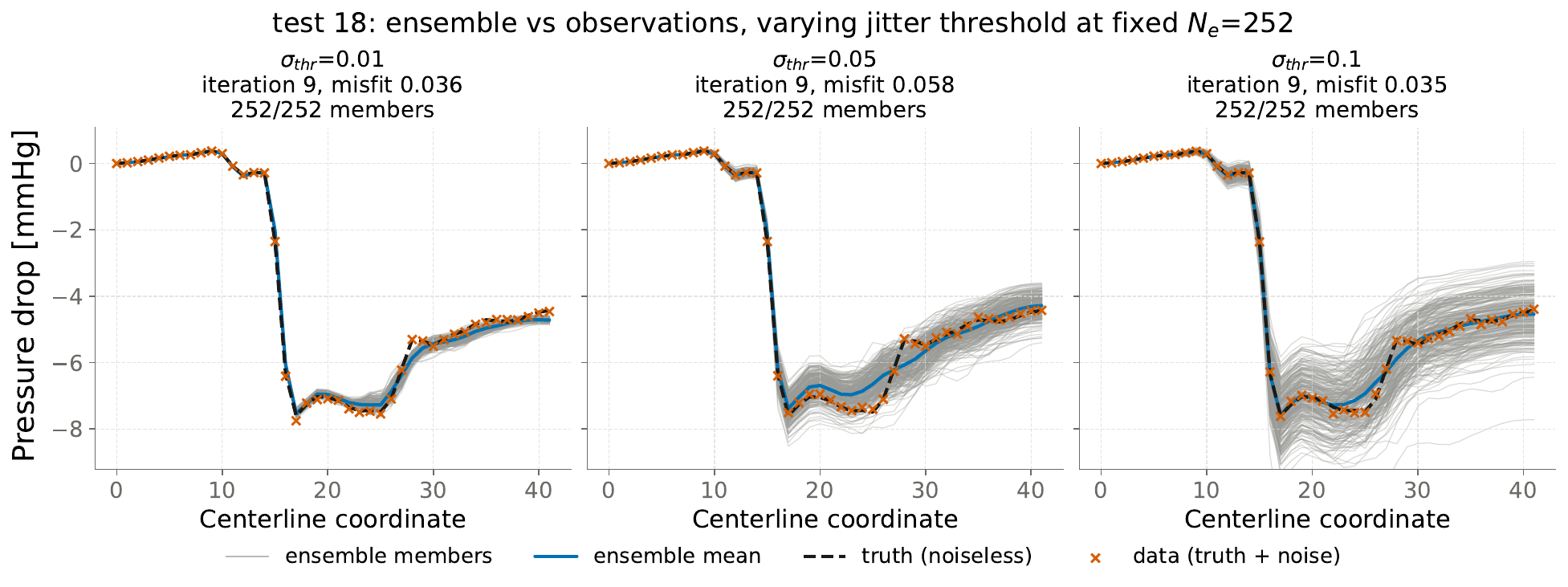}
  \caption{Test case $18$. Centerline pressure drops for the EAKI ensemble at the final iteration ($T_{\text{max}}=10$), varying the jittering threshold $\sigma_{\text{thr}}$.}
  \label{fig:eakithr}
\end{figure}

\begin{figure}
  \includegraphics[width=\textwidth]{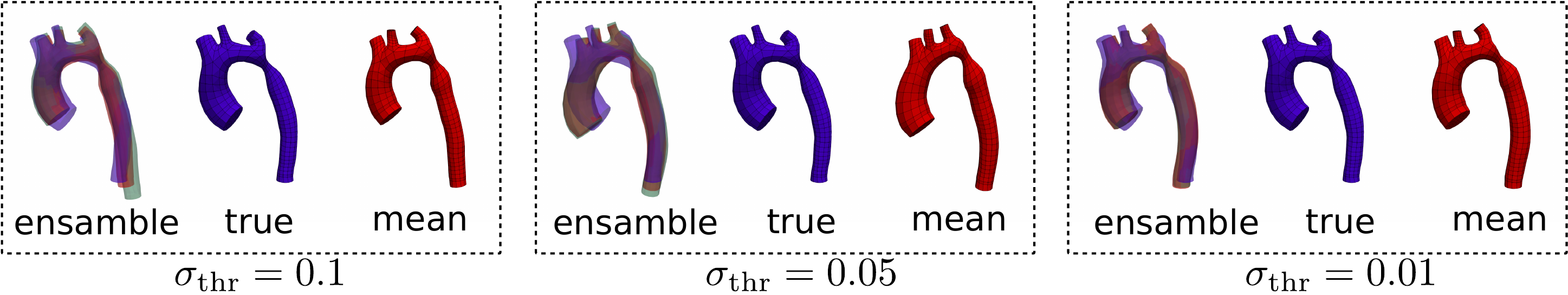}
  \caption{Test case $18$. Jittering results corresponding to figure~\ref{fig:eakithr}. For each value of 
  $\sigma_{\text{thr}}$, the picture shows the true (blue), the mean (red), and the farthest (green) EAKI ensemble members overlaid (left), the true geometry (center), and the EAKI mean geometry (right). The ensemble collapses to its mean from $\sigma_{\text{thr}}=0.1$ to $\sigma_{\text{thr}}=0.01$.}
  \label{fig:eakicollapse}
\end{figure}

The results of jittering with $\sigma_{\text{thr}}=0.1, 0.05, 0.01$ are shown in figures~\ref{fig:eakithr} and~\ref{fig:eakicollapse}. The jittering can be used to avoid ensemble collapse, and one can see that the ensemble variance reduces considerably when decreasing the threshold from $\sigma_{\text{thr}}=0.1$ to $\sigma_{\text{thr}}=0.01$. A higher ensemble variance is useful in practice to better study or model the uncertainty of the shape inversion problem: since we are evaluating the pressure drop with respect to the inlet pressure, the uncertainty is mostly distributed in the descending aortic arch rather than in the ascending arch, which is closer to the inlet, as can be seen in figure~\ref{fig:eakicollapse}.

\begin{figure}
  \includegraphics[width=\textwidth]{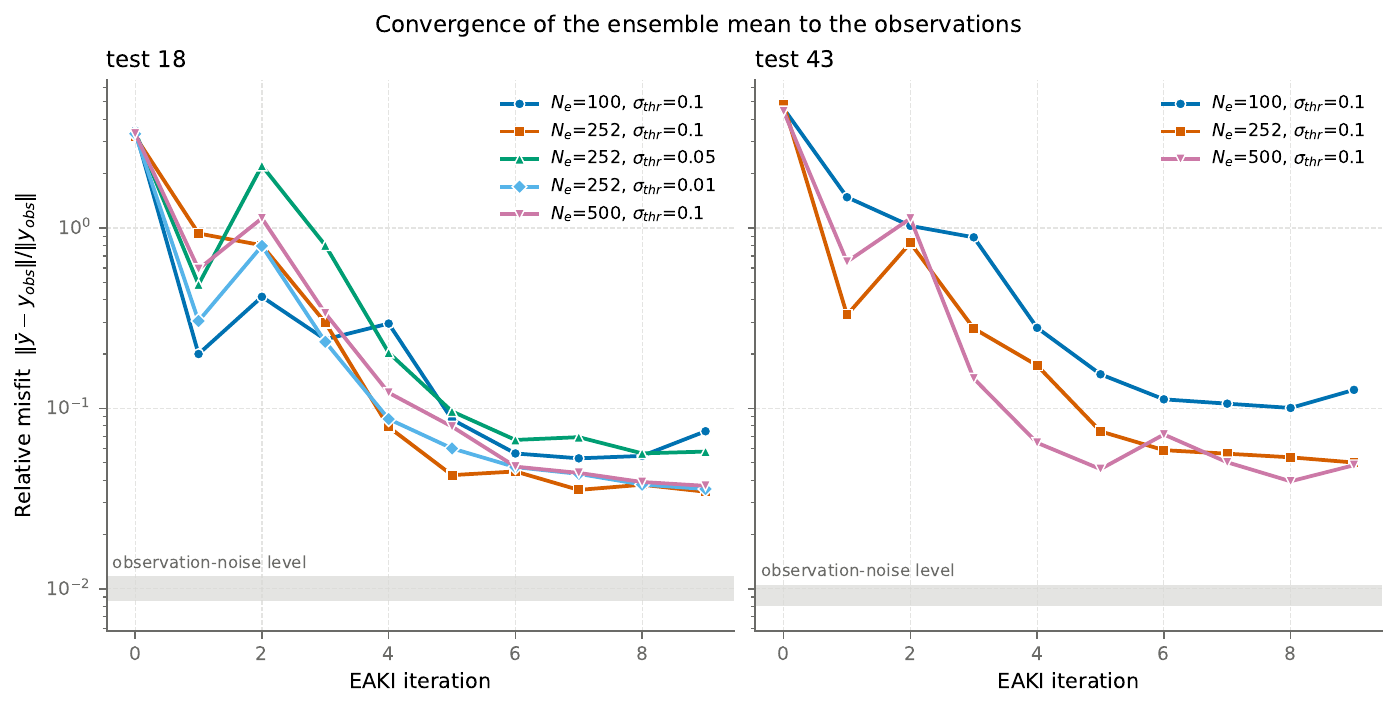}
  \caption{Misfit $\lVert y-\mathcal{G}(\theta)\rVert_{2}/\lVert y\rVert_{2}$ over
  the iterations. Left: test case $18$. Right: test case $43$.}
  \label{fig:eaki}
\end{figure}

The convergence of EAKI with respect to the relative misfit $\lVert y-\mathcal{G}(\theta)\rVert_{2}/\lVert y\rVert_{2}$ is shown for the considered cases in figure~\ref{fig:eaki}. It can be seen that the pressure drop profile at time $t=0.25\,\mathrm{s}$ is sufficient to identify a shape close to the true one after only $6$ iterations. 

The computational costs for the EAKI are reported in table~\ref{tab:uq-cost}, detailing the ensemble size, the number of evaluations of the forward map, the wall-clock time (total and per evaluation), and the CPU time (in $32$-core node hours).

A few forward evaluations did not terminate correctly in the two runs with $N_{\mathrm{ens}}=500$ (only $4963$ and $4955$ valid evaluations instead of $5000$). In fact, on a few meshes with particularly pronounced geometrical features, which require a very small adaptive time step, the numerical solver did not converge. These failures occurred only at the first EAKI iteration, since only in that case the particles are generated by Latin hypercube sampling, and not for the other considered methods.

\subsection{Scalar pressure drop}
\label{subsec:scalar_pressure_drop}
In the next test, we consider the shape inversion problem based on a single observed value corresponding 
to the scalar pressure drop $\mathcal{O}_2$: the values are chosen by hand from $\{5\,\mathrm{mmHg}, 10\,\mathrm{mmHg}, 20\,\mathrm{mmHg}\}$ and corresponding to cases ranging from healthy to severe conditions. For this setting, we study the performance of the aBUS and the WFR methods. %

\subsubsection{aBUS}
\label{subsec:results_aortic_abus}
The main parameter that affects the convergence of aBUS is the number of samples ($N$). For the other parameters, 
we use the values suggested in the literature~\cite{simon2026parallelized}, fixing
the target acceptance rate $\hat{a}_{\text{target}}=0.44$, the conditional probability between consecutive levels $p_0=0.1$, and the batch size $n_B=100$.

\begin{figure}[htp!]
  \includegraphics[width=\textwidth]{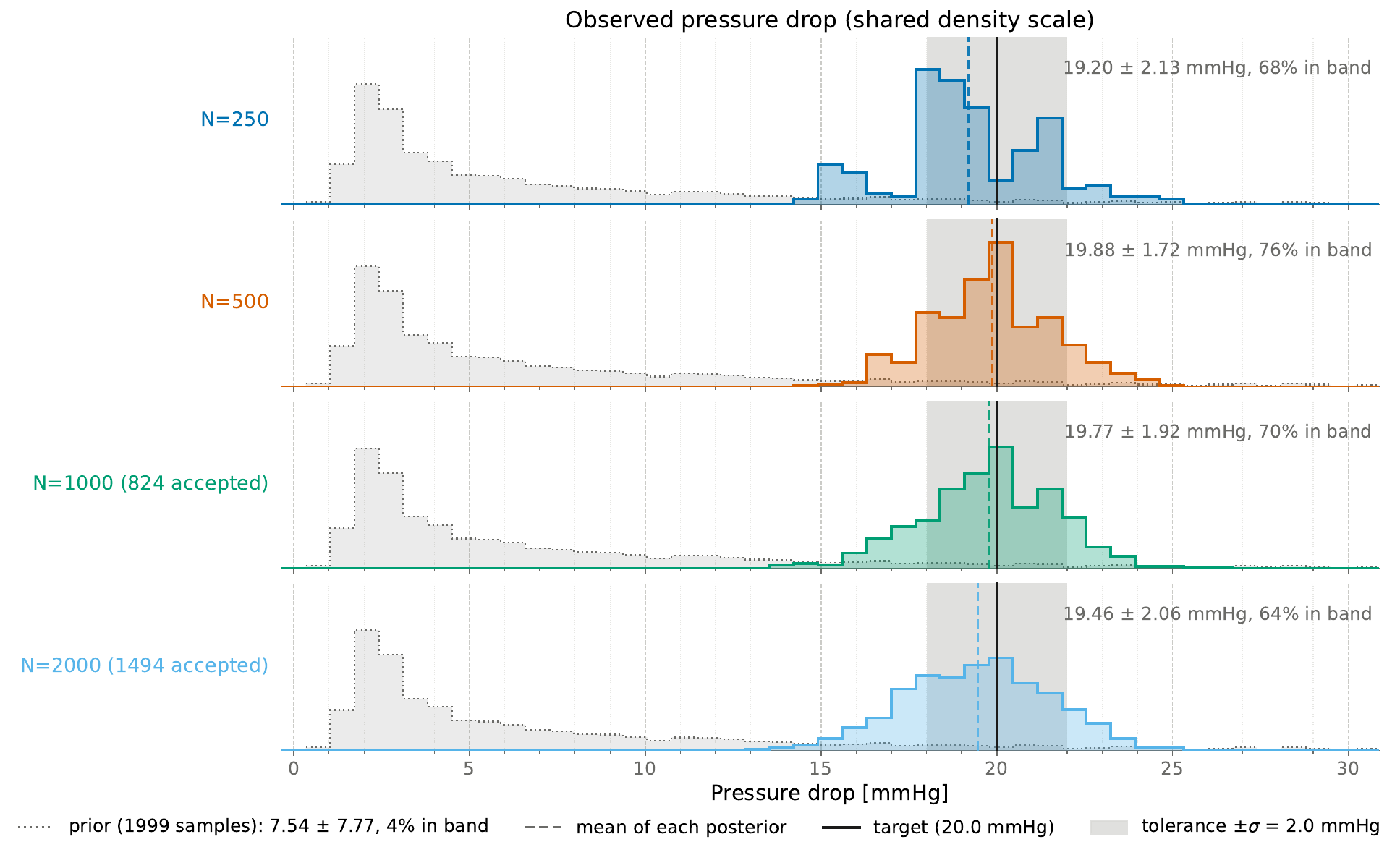}
  \caption{Convergence of the aBUS method with respect to the sample size $N=250, 500, 1000, 2000$, with $\mathcal{O}_2=20\,\mathrm{mmHg}$.}
  \label{fig:convergence_abus}
\end{figure}

We consider a pressure drop of $20\,\mathrm{mmHg}$, for which the $10\%$ observation noise gives $\sigma_\eta=2\,\mathrm{mmHg}$.
The convergence with respect to the number of samples $N=250, 500, 1000, 2000$ is shown in figure~\ref{fig:convergence_abus}.
This pressure drop is considerably larger than physiological values for healthy subjects, and 
it can be obtained only in the presence of an additional stenosis. Figure~\ref{fig:convergence_abus} shows also the values of the observations associated to samples from the prior %
distribution (figure~\ref{fig:surrogate_dist}).

\begin{figure}[htp!]
  \includegraphics[width=0.8\textwidth]{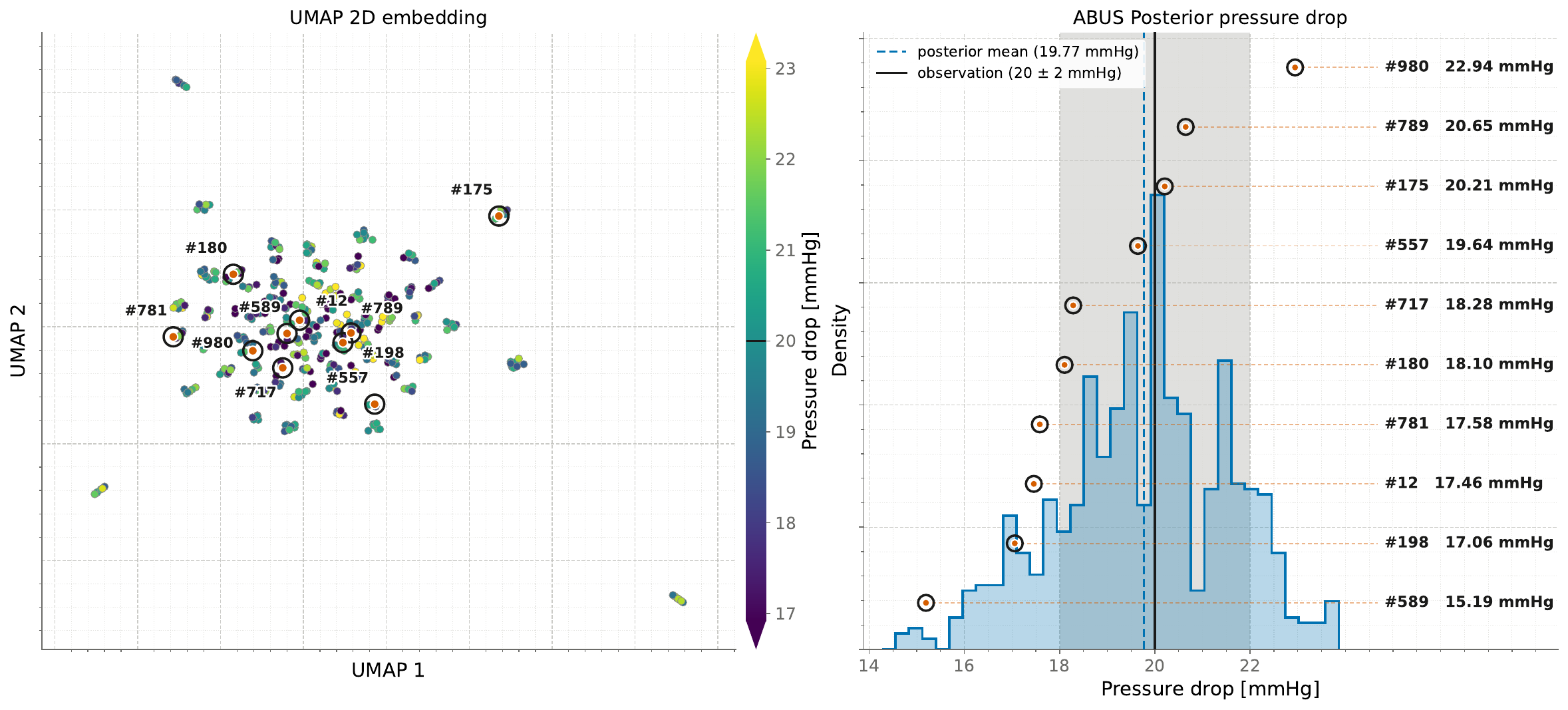}\\
   \includegraphics[width=1\textwidth]{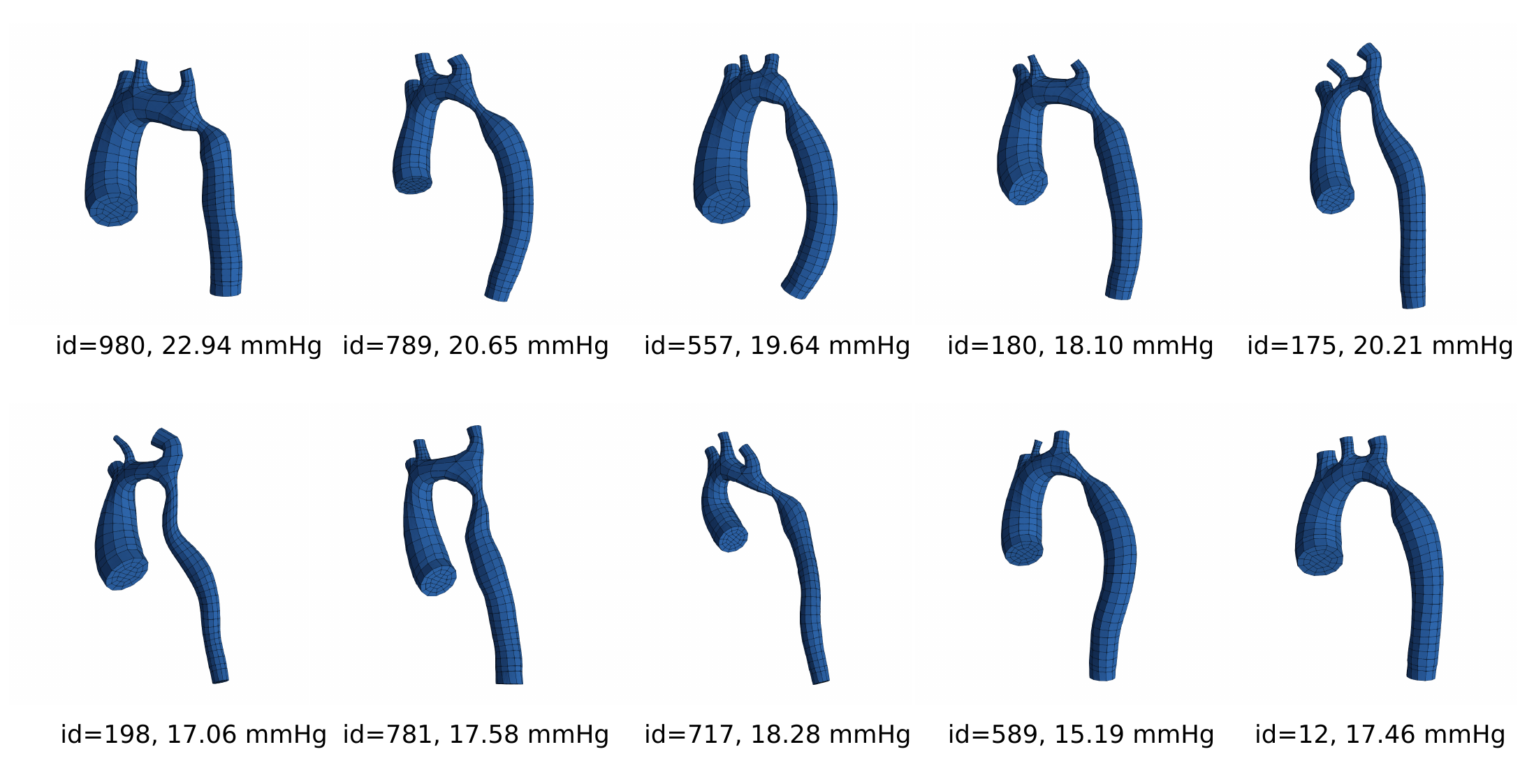}
  \caption{aBUS shape inversion with maximum pressure drop $\mathcal{O}_2=20\,\mathrm{mmHg}$. Posterior distribution 
  from aBUS and selected samples.}
  \label{fig:multimodal_abus}
\end{figure}

The posterior distribution, as well as
selected samples are shown in figure~\ref{fig:multimodal_abus}. In particular, figure~\ref{fig:multimodal_abus}, top-left, depicts the samples embedded in $\mathbb{R}^2$ with Uniform Manifold Approximation and Projection (UMAP).
Unlike EAKI, the samples from the aBUS posterior distribution comprise multiple qualitatively different geometries rather than a single ensemble-averaged one.

\subsubsection{WFR}
\label{subsec:results_aortic_wfr}
The WFR uses the surrogate forward map (section~\ref{sec:forward_lddmm}), which does not include the additional stenosis. The method is
therefore tested considering two moderate values of pressure drop, $\mathcal{O}_2\in\{5,10\}\,\mathrm{mmHg}$, which do not lie too far from the prior/training distribution on which the surrogate model is trained (figure~\ref{fig:surrogate_dist}). Without the stenosis, higher values of the pressure drop 
would correspond to a worse accuracy of the surrogate and, therefore, of the posterior distribution.

The main parameter considered for the convergence study is the number of samples $N=250, 500, 1000$. Other relevant parameters are the number of steps of the gradient flow dynamics (fixed to $T_{\text{max}}=10$), the minimum and maximum values of the cosine scheduler (fixed as $\gamma_\text{min}=0.01$, $\gamma_\text{max}=0.1$), and the threshold on the effective sample size (ESS) set to $N/2$.

\begin{figure}[htp!]
  \includegraphics[width=\textwidth]{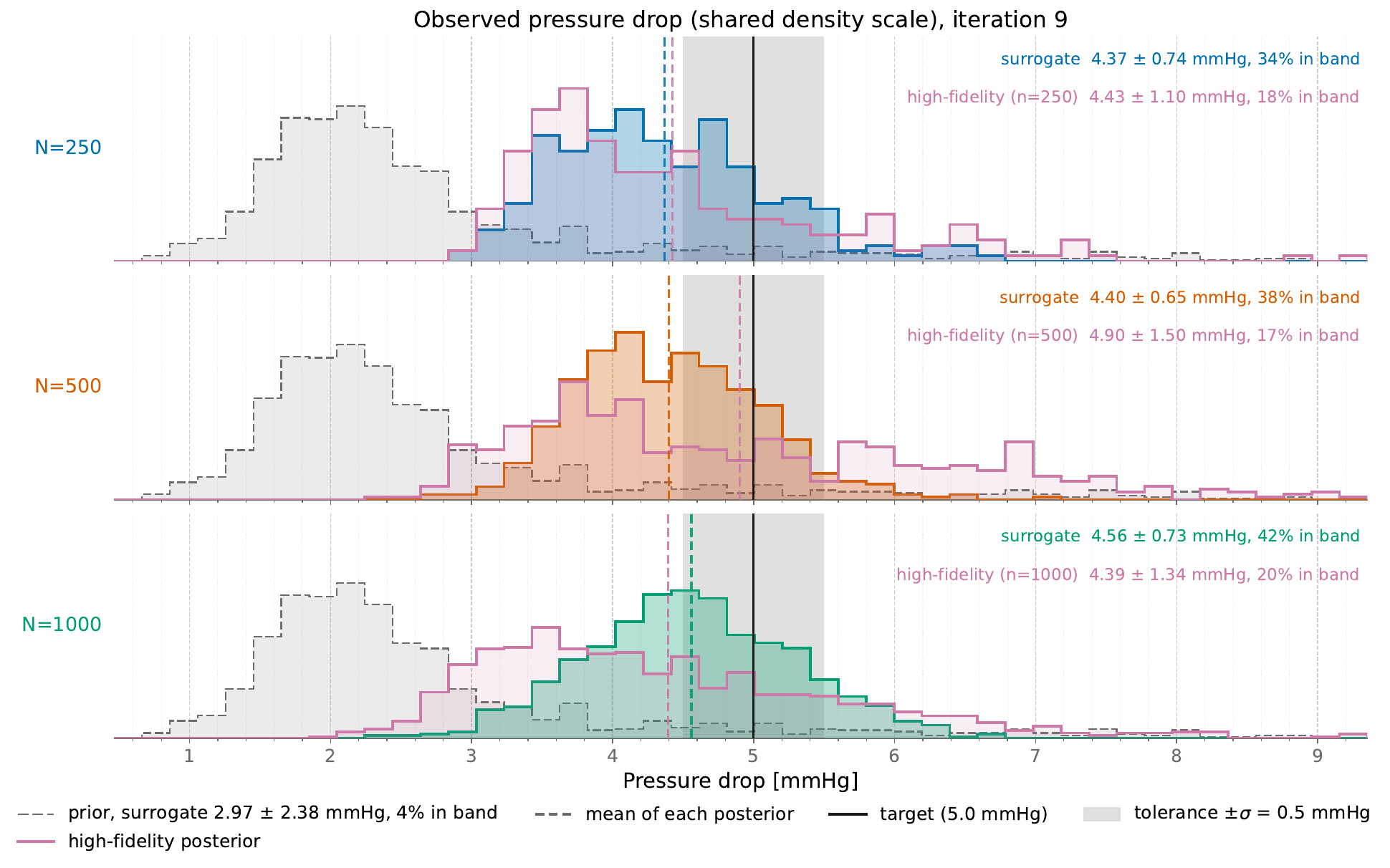}
  \caption{WFR, scalar pressure drop $\mathcal{O}_2=5\,\mathrm{mmHg}$. Convergence of the WFR posterior distribution with respect to the ensemble size $N=250, 500, 1000$ (from top to bottom).}
  \label{fig:convergence_wfr_5}
\end{figure}

\begin{figure}[htp!]
  \includegraphics[width=\textwidth]{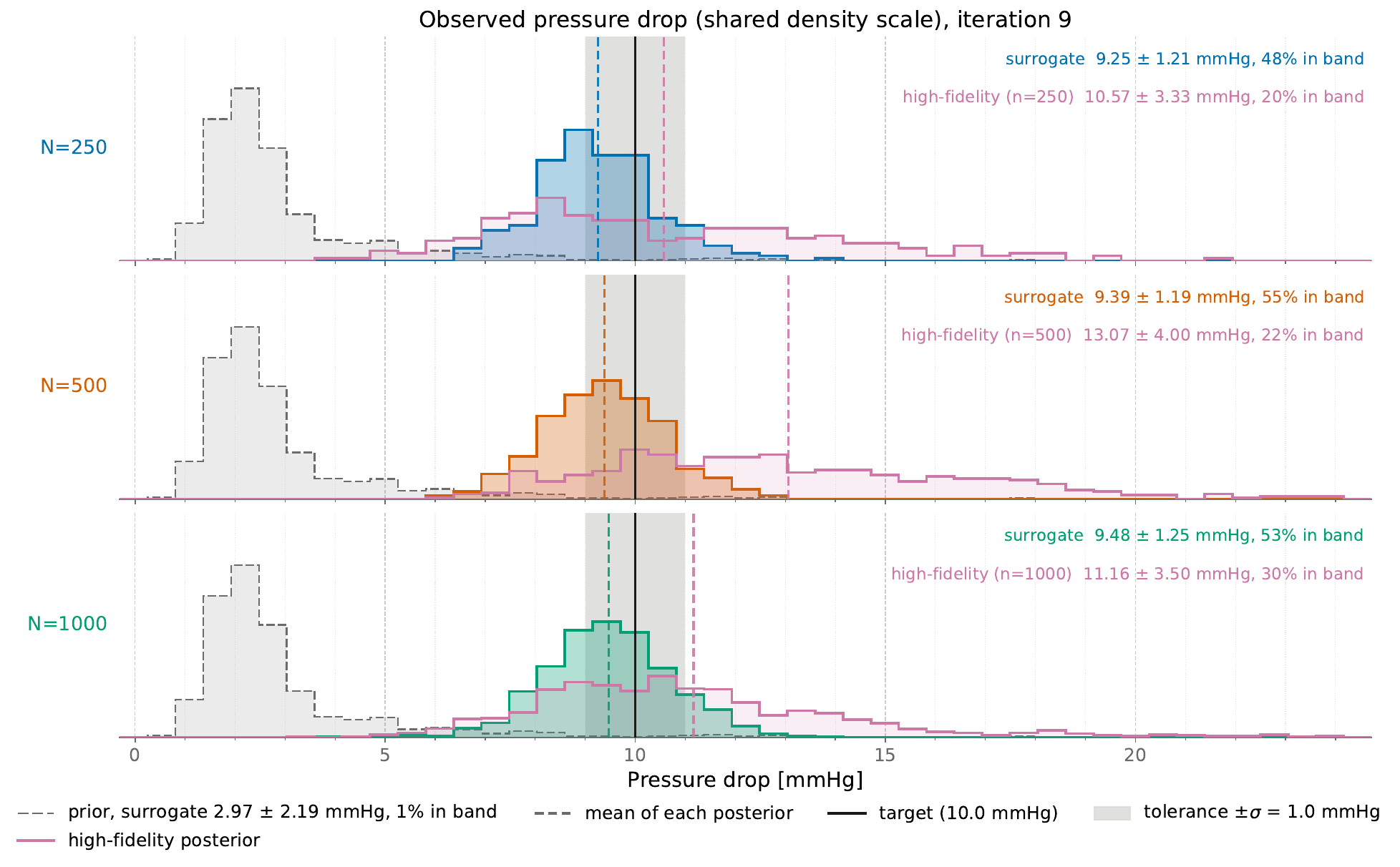}
  \caption{WFR, scalar pressure drop $\mathcal{O}_2=10\,\mathrm{mmHg}$. Convergence of the WFR posterior distribution with respect to the ensemble size $N=250, 500, 1000$ (from top to bottom).}\label{fig:convergence_wfr_10}
\end{figure}

Figures~\ref{fig:convergence_wfr_5} and~\ref{fig:convergence_wfr_10} show the 
results for pressure drops of $5\,\mathrm{mmHg}$ and $10\,\mathrm{mmHg}$, respectively.
After sampling the parameters $\Tilde{\theta}_j\in\mathbb{R}^{n_c}$ from the WFR posterior, the pressure values are computed by generating the corresponding computational domain and computing the pressure with the high-fidelity solver.

Notice that values of the pressure drop farther from those of the training distribution correspond to a lower accuracy of the surrogate. As a result, when calculating the actual pressure drop corresponding to the posterior distribution, we observe a distribution shift. The case $\mathcal{O}_2=5\,\mathrm{mmHg}$ is less affected than $\mathcal{O}_2=10\,\mathrm{mmHg}$. 
As observed in section~\ref{sec:forward_lddmm} (see figure~\ref{fig:rbf_sensitivity}), this shift
is not due to the 
different approach used to generate the volume mesh, but should rather be attributed to the sensitivity of the surrogate to the training dataset.


These results can be improved by increasing the amount of training data associated with the prescribed values of the pressure drop, or by preventing the WFR method from exploring regions of the parameter space too far from the prior distribution. 
Nonetheless, the predicted posterior distributions match the prescribed values with good accuracy, which indicates that the distribution shift observed with the high-fidelity values is mainly due to the low accuracy of the surrogate model.

\begin{figure}[htp!]
  \includegraphics[width=.8\textwidth]{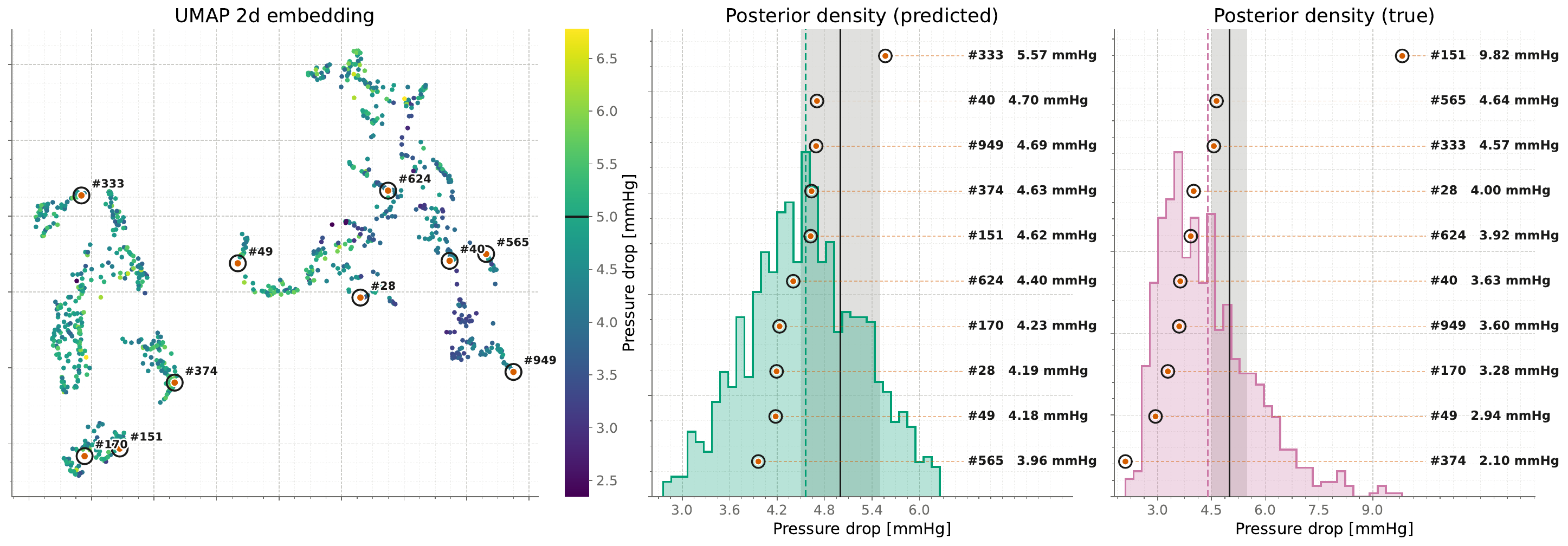}\\
  \includegraphics[width=.8\textwidth]{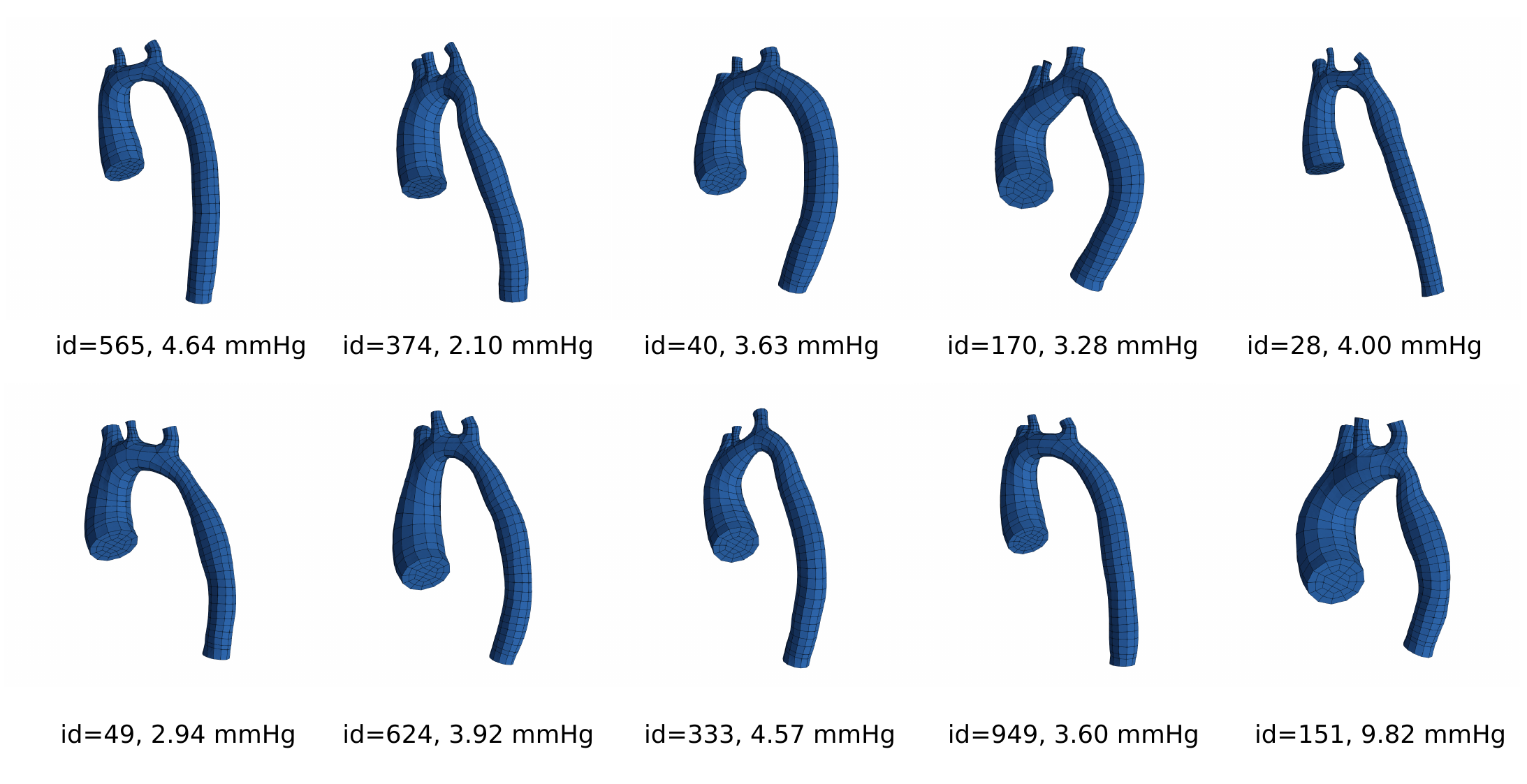}
  \caption{Representative sampled geometries from the WFR posterior distribution with $\mathcal{O}_2=5\,\mathrm{mmHg}$.}
  \label{fig:cluster_wfr_5}
\end{figure}

\begin{figure}[htp!]
  \includegraphics[width=.8\textwidth]{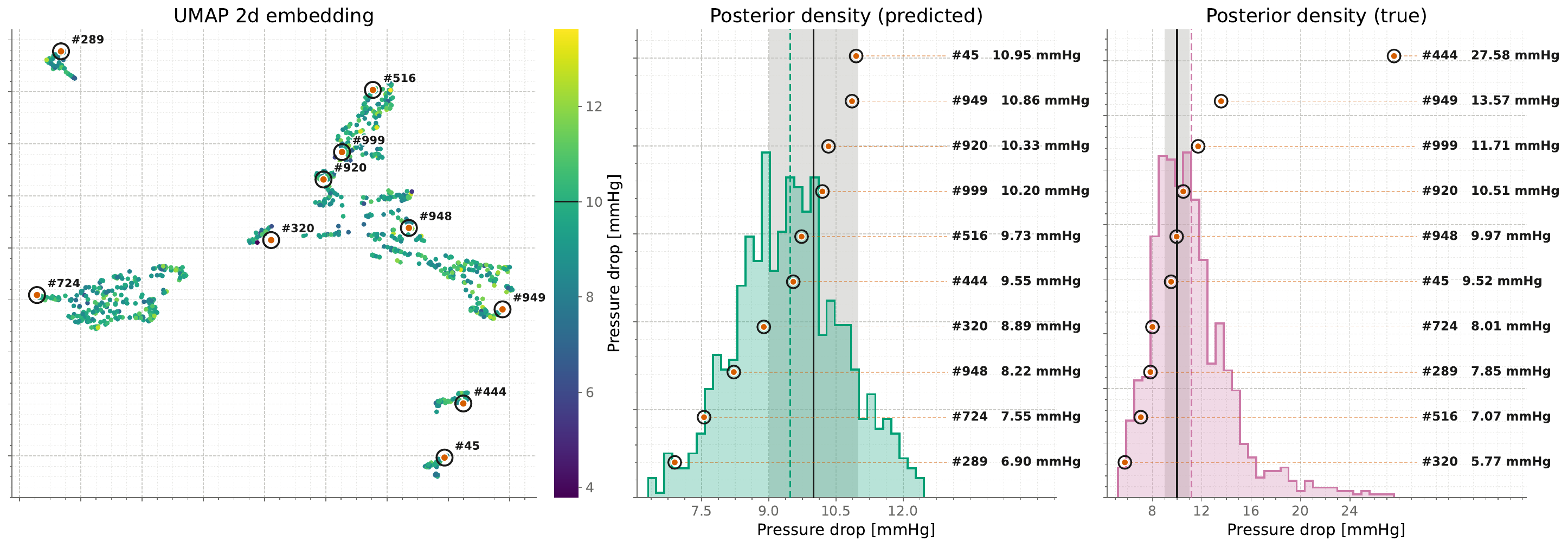}\\
  \includegraphics[width=.8\textwidth]{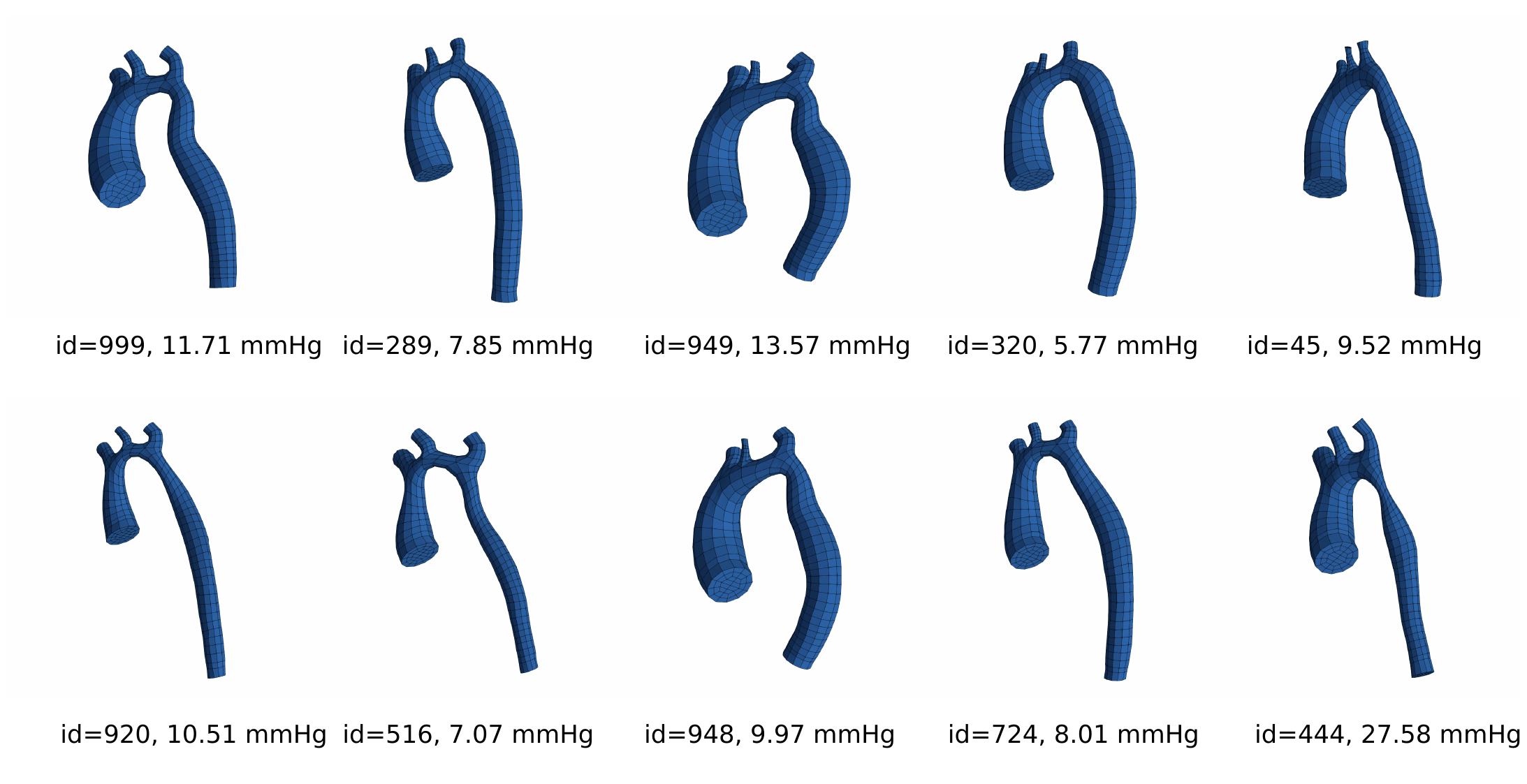}
  \caption{Representative sampled geometries from the WFR posterior distribution with pressure drop of $10\,\mathrm{mmHg}$.}
  \label{fig:cluster_wfr_10}
\end{figure}

As for aBUS, the WFR posterior distributions
can be visualized in $\mathbb{R}^2$ with
a UMAP embedding, see figures~\ref{fig:cluster_wfr_5} and~\ref{fig:cluster_wfr_10}.

\subsubsection{Comparisons: aBUS vs. WFR}
\begin{figure}[htp!]
\centering
  \includegraphics[width=.9\textwidth]{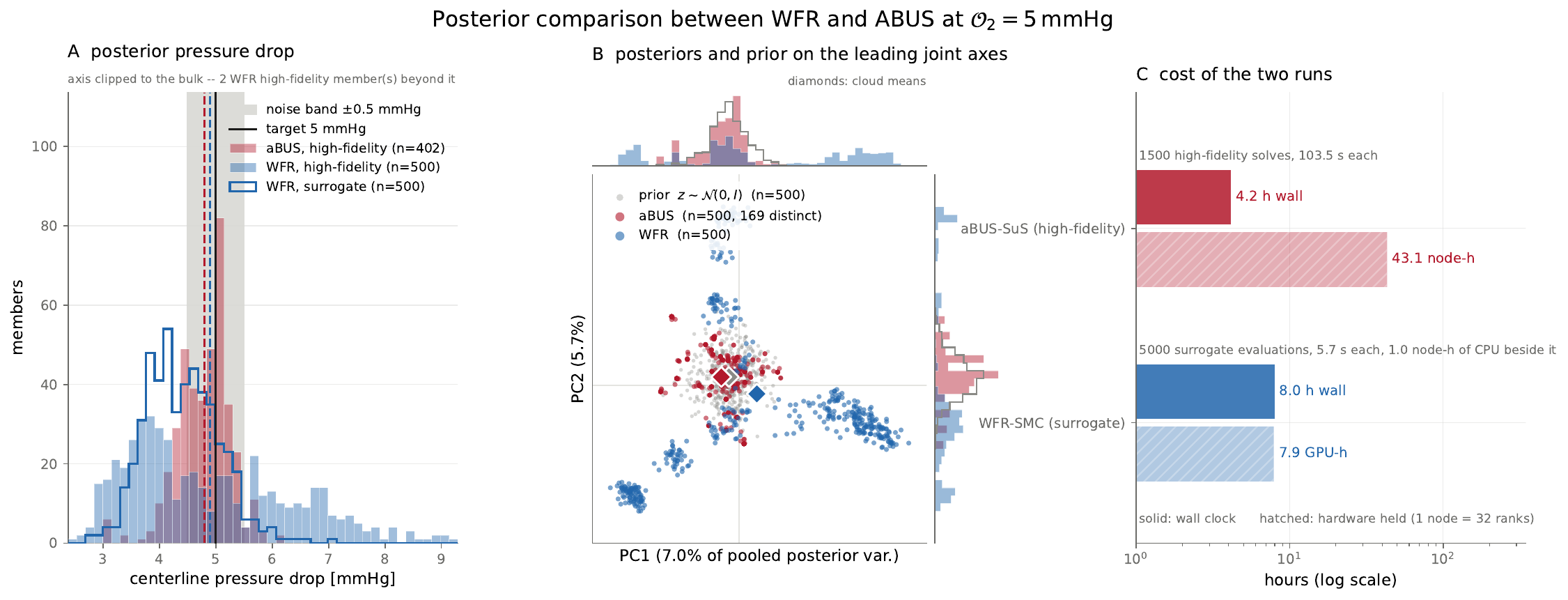}\\
  \includegraphics[width=.9\textwidth]{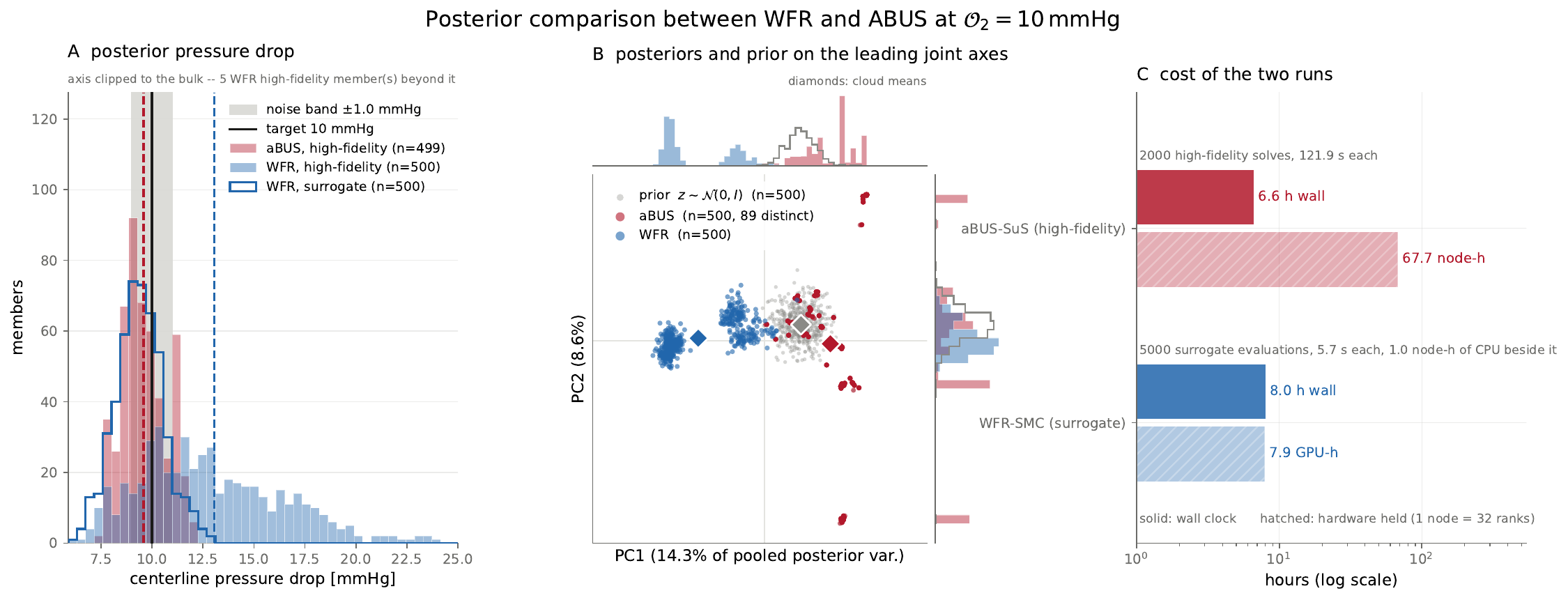}\\
  \caption{Comparison of the posterior distributions of WFR and aBUS for pressure drops of $5\,\mathrm{mmHg}$ (top) and $10\,\mathrm{mmHg}$ (bottom).}
  \label{fig:compare}
\end{figure}

We can compare the aBUS and WFR posterior distributions by restricting to the parameter space
without the stenosis ($\alpha=0$).
The resulting posterior distributions for pressure drops of $5\,\mathrm{mmHg}$ and $10\,\mathrm{mmHg}$ are shown in figure~\ref{fig:compare}, considering the same number of samples ($N=500$).
In terms of misfit, aBUS is more accurate. WFR is accurate on the predicted posterior values but not on the posterior values computed with the high-fidelity solver. Hence, the higher misfit of the WFR results is mainly due to the limited accuracy of the surrogate outside the training distribution. This can be seen in the central plots in figure~\ref{fig:compare}, where the aBUS, the WFR, and the prior distributions are shown embedded in the first two principal components of a joint PCA, showing that aBUS remains much closer to the prior distribution in the parameter space.

In terms of computational costs, 
table~\ref{tab:uq-cost} provides a detailed
comparison of the methods, considering the ensemble size, the number
of evaluations of the forward map, the wall-clock time (total and per evaluation), and the CPU time (in $32$-core node hours).

Using the surrogate forward map considerably reduces the computational cost. The HTCondor queue has $4$ x64/openSUSE15 nodes with $72, 128, 64, 128$ cores, respectively: each high-fidelity solver run requires $32$ CPU cores, and a single CPU node can run multiple forward maps simultaneously. The surrogate runs in batches on a single NVIDIA H200 NVL4 141GB GPU. 

Although this imbalance of resources results in a shorter wall-clock time for aBUS ($\sim 4$--$6$ hours against $8$ hours for WFR), the compute time, measured in $32$-core node hours and GPU hours, is higher for aBUS than for WFR, with an average forward evaluation time of $5.7\,\mathrm{s}$ for WFR against around $120\,\mathrm{s}$ for aBUS. 
Notice that the WFR steps are performed on CPU, while only the forward and gradient evaluations of the surrogate are performed on GPU.
For example, compare the 
$67.7$ node hours required by aBUS for the setting
with pressure drop $10\,\mathrm{mmHg}$ and ensemble size $N_{\mathrm{ens}}=500$ ($2000$ forward map evaluations), against $0.5$ node hours plus $4$ GPU hours spent by WFR for pressure drop $10\,\mathrm{mmHg}$ with $N_{\mathrm{ens}}=250$, which requires more forward map evaluations ($2500$), including also the gradients of the observation operator.

\begin{table}[htp!]
  \centering
  \caption{Computational cost of the aortic Bayesian inversion studies. $N_{\mathrm{ens}}$ is the ensemble size
    ($N$ for EAKI, samples per level for aBUS, number of particles for WFR), $N_{\mathrm{fwd}}$ the total number of
    forward map evaluations, $T_{\mathrm{wall}}$ the wall-clock time of the whole study, $\bar t_{\mathrm{fwd}}$ the
    mean wall-clock time of a single forward map evaluation, and $T_{\mathrm{wall}}/N_{\mathrm{fwd}}$ the same time
    amortized over the study, which measures how effectively the evaluations are executed concurrently. The last two
    columns report the consumed resources: CPU time normalized to $32$-core nodes and GPU time. 
    For EAKI, a few forward evaluations did not terminate correctly, so that $N_{\text{fwd}}$ is lower than $10\,N_{\mathrm{ens}}$.
}
  \label{tab:uq-cost}
  \small
  \setlength{\tabcolsep}{5pt}
  \begin{tabular}{lrrrrrrr}
    \toprule
    & $N_{\mathrm{ens}}$ & $N_{\mathrm{fwd}}$ & $T_{\mathrm{wall}}$ & $\bar t_{\mathrm{fwd}}$ & $T_{\mathrm{wall}}/N_{\mathrm{fwd}}$ & CPU & GPU \\
    Run & & & [h:mm] & [s] & [s] & [$32$-core node h] & [h] \\
    \midrule
    \multicolumn{8}{@{}l}{\textbf{EAKI}\quad\textit{$\mathcal{O}_1$, $T_{\text{max}}=10$ iterations, one high-fidelity solve per ensemble member}} \\
    \addlinespace[1pt]
    \quad test $18$, $\sigma_{\text{thr}}=0.1$ & $100$ & $1000$ & $3{:}42$ & $114.0$ & $13.3$ & $31.7$ & --- \\
    \quad test $18$, $\sigma_{\text{thr}}=0.01$ & $252$ & $2520$ & $6{:}55$ & $115.8$ & $9.9$ & $81.1$ & --- \\
    \quad test $18$, $\sigma_{\text{thr}}=0.05$ & $252$ & $2520$ & $14{:}29$ & $113.1$ & $20.7$ & $79.2$ & --- \\
    \quad test $18$, $\sigma_{\text{thr}}=0.1$ & $252$ & $2520$ & $7{:}20$ & $118.4$ & $10.5$ & $82.9$ & --- \\
    \quad test $18$, $\sigma_{\text{thr}}=0.1$ & $500$ & $4963$ & $13{:}23$ & $120.2$ & $9.7$ & $166.2$ & --- \\
    \quad test $43$, $\sigma_{\text{thr}}=0.1$ & $100$ & $1000$ & $2{:}52$ & $110.1$ & $10.3$ & $30.6$ & --- \\
    \quad test $43$, $\sigma_{\text{thr}}=0.1$ & $252$ & $2520$ & $6{:}35$ & $111.8$ & $9.4$ & $78.2$ & --- \\
    \quad test $43$, $\sigma_{\text{thr}}=0.1$ & $500$ & $4955$ & $11{:}34$ & $103.7$ & $8.4$ & $143.4$ & --- \\
    \midrule
    \multicolumn{8}{@{}l}{\textbf{aBUS--SuS}\quad\textit{$\mathcal{O}_2$, one high-fidelity solve per sample and level}} \\
    \addlinespace[1pt]
    \quad $\mathcal{O}_2=5\,\mathrm{mmHg}$ & $500$ & $1500$ & $4{:}09$ & $103.5$ & $10.0$ & $43.1$ & --- \\
    \quad $\mathcal{O}_2=10\,\mathrm{mmHg}$ & $500$ & $2000$ & $6{:}36$ & $121.9$ & $11.9$ & $67.7$ & --- \\
    \quad $\mathcal{O}_2=20\,\mathrm{mmHg}$ & $250$ & $1000$ & $4{:}03$ & $147.6$ & $14.6$ & $41.0$ & --- \\
    \quad $\mathcal{O}_2=20\,\mathrm{mmHg}$ & $500$ & $2000$ & $16{:}21$ & $155.4$ & $29.4$ & $86.3$ & --- \\
    \quad $\mathcal{O}_2=20\,\mathrm{mmHg}$ & $1000$ & $3000$ & $10{:}08$ & $129.3$ & $12.2$ & $107.8$ & --- \\
    \quad $\mathcal{O}_2=20\,\mathrm{mmHg}$ & $2000$ & $6000$ & $18{:}18$ & $117.7$ & $11.0$ & $196.1$ & --- \\
    \midrule
    \multicolumn{8}{@{}l}{\textbf{WFR}\quad\textit{$\mathcal{O}_2$, $T_{\max}=10$ iterations, one batched GPU job per iteration}} \\
    \addlinespace[1pt]
    \quad $\mathcal{O}_2=5\,\mathrm{mmHg}$ & $250$ & $2500$ & $4{:}03$ & $5.7$ & $5.8$ & $0.5$ & $4.0$ \\
    \quad $\mathcal{O}_2=5\,\mathrm{mmHg}$ & $500$ & $5000$ & $7{:}57$ & $5.7$ & $5.7$ & $1.0$ & $7.9$ \\
    \quad $\mathcal{O}_2=5\,\mathrm{mmHg}$ & $1000$ & $10000$ & $15{:}43$ & $5.6$ & $5.7$ & $2.0$ & $15.6$ \\
    \quad $\mathcal{O}_2=10\,\mathrm{mmHg}$ & $250$ & $2500$ & $4{:}03$ & $5.7$ & $5.8$ & $0.5$ & $4.0$ \\
    \quad $\mathcal{O}_2=10\,\mathrm{mmHg}$ & $500$ & $5000$ & $7{:}58$ & $5.7$ & $5.7$ & $1.0$ & $7.9$ \\
    \quad $\mathcal{O}_2=10\,\mathrm{mmHg}$ & $1000$ & $10000$ & $15{:}50$ & $5.7$ & $5.7$ & $2.0$ & $15.7$ \\
    \bottomrule
  \end{tabular}
\end{table}

\section{Conclusions}\label{sec:conclusions}
We addressed an inverse domain uncertainty quantification problem to determine virtual patient cohorts from target biomarkers in aortic coarctation.
This approach has potential applications for data augmentation in the context of in silico testing of medical devices and treatments on patients with prescribed pathological conditions.
Moreover, framed as a shape optimization problem, it can be 
applied to the estimation of the effect of geometric uncertainty introduced by the acquisition process of MRI and CT scans on numerical simulations.
Using an initial cohort of synthetic patients, we define a latent
parameter space ($250$-dimensional) 
that implicitly encodes centerlines and inscribed radii of the given population. The forward map for the generation of shapes is based on a ResNet-LDDMM
aortic shape registration method~\cite{romor2026shape}, combined with
LDDMM stochastic interpolants generative modeling~\cite{KATZ2026119120}, trained on the available set of shapes and corresponding pairwise registration maps.
We considered three of the most promising UQ methods for 
large-scale and high-dimensional applications:
EAKI, which relies on Gaussian approximations of the sought
posterior distribution, and aBUS and WFR, which sample from multimodal and high-dimensional posterior distributions. 
We studied two different scenarios, considering observations of the pressure drop along the centerline (EAKI) or only its maximum value (aBUS and WFR).
The results show that the posterior distributions reproduce the
target measurements in all considered cases.
Moreover, we observed that EAKI is more suitable for quickly extracting a single geometry with prescribed biomarkers from the generative model, whereas aBUS and WFR can be used to sample an entire cohort of aortic shapes associated with a prescribed biomarker of interest.

The present studies were conducted considering only the pressure drop. However, the general framework can be extended to other quantities, such as the velocity field and the wall shear stress, among others.
The proposed methodology is also not limited to the considered application:  it can be applied to any family of topologically equivalent shapes for which a training dataset is available. In this setting, registration maps onto a common template provide a natural geometrical encoding, shared by the generative prior and the surrogate model, and remove the need for a hand-crafted shape parametrization.

The restriction to a single topology is, at the same time, one of the main limitations of the framework: extending the generative prior to shapes with different topologies, e.g., aortas with a varying number of supra-aortic branches, requires different geometrical encodings.
A further limitation concerns geometrical constraints that the generative prior does not enforce by construction, such as the absence of self-intersections or a minimum vessel radius. These constraints are currently satisfied only to the extent that they are represented in the training data, while invalid geometries are handled a posteriori by the forward map (remark~\ref{rmk:invalid}). Incorporating them explicitly, for instance through constrained generative models or additional terms in the likelihood, is a further direction of future research.

From a broader perspective, this work fits into the emerging field of generative design, also referred to as AI-assisted design, in which a generative model of admissible shapes is coupled with a physics-based forward model to find, or sample, geometries that attain prescribed performance targets.
Further  research directions 
include the improvement of the Fisher--Rao step of the WFR method through a more accurate particle density estimator, and the validation of the surrogate model with more training samples, which are needed to improve its accuracy.

\section*{Acknowledgements}
Funded by the Deutsche Forschungsgemeinschaft (DFG, German Research
Foundation) under Germany's Excellence Strategy – The Berlin Mathematics
Research Center MATH+ (EXC-2046/1, EXC-2046/2, project ID: 390685689).

The authors used Claude Code (Anthropic) to polish the text and to check its grammar, and in part to support the implementation of the uncertainty quantification methods. The authors have reviewed and edited the output and take full responsibility for the content of this work.

\printbibliography
\newpage
\appendix

\end{document}